\documentclass[english,11pt,twoside]{article}
\usepackage{babel,a4}
\usepackage{amsmath,amssymb}
\usepackage[latin1]{inputenc}
\usepackage{graphicx,rotate}
\usepackage{palatino}
\usepackage{array}
\usepackage{latexsym}
\usepackage{mathbbol}
\numberwithin{equation}{section}
\title{A Weak Limit Shape Theorem For Planar Isotropic Brownian Flows}
\author{Holger van Bargen}
\date{Version: \today}
\newcommand{\ind}[1]{\ensuremath{ \mathbb{1}_{{#1}} } }
\newcommand{\indb}[1]{\ensuremath{ \mathbb{1}_{\left\{ {#1}\right\} } } }

\newcommand{\ce}[3][]{\ensuremath{\mathbb{E}_{{#1}}\left[{#2}\left|{#3}\right.\right]}}
\newcommand{\cel}[3][]{\ensuremath{\mathbb{E}_{{#1}}\left[\left.{#2}\right|{#3}\right]}}
\newcommand{\cp}[3][]{\ensuremath{\mathbb{P}_{{#1}}\left[{#2}\left|{#3}\right.\right]}}

\newcommand{\dist}{\textnormal{dist}}
\newcommand{\diam}{\textnormal{diam}}
\newcommand{\be}{\begin{equation}}
\newcommand{\ee}{\end{equation}}
\newcommand{\bee}{\begin{eqnarray}}
\newcommand{\eee}{\end{eqnarray}}
\newcommand{\ip}[2]{\ensuremath{\left\langle{#1},{#2}\right\rangle}}
\newcommand{\prob}[2][]{\ensuremath{\mathbb{P}_{{#1}}\left[{#2}\right]}}
\newcommand{\expec}[2][]{\ensuremath{\mathbb{E}_{{#1}}\left[{#2}\right]}}
\newtheorem{theorem}{Theorem}[section]
\newtheorem{lemma}[theorem]{Lemma}
\newtheorem{definition}[theorem]{Definition}

\newtheorem{proposition}[theorem]{Proposition}
\newtheorem{corollary}[theorem]{Corollary}

\newtheorem{condition}[theorem]{Condition}

\newcommand{\R}{\ensuremath{\mathbb{R}}}
\newcommand{\N}{\ensuremath{\mathbb{N}}}
\newcommand{\Rd}{\ensuremath{\mathbb{R}^{d}}}
\newcommand{\eps}{\ensuremath{\epsilon}}
\begin{document}
\thispagestyle{empty}
\vglue20pt \centerline{\Large\bf A Weak Limit Shape Theorem For Planar Isotropic Brownian Flows }

\medskip

\bigskip

\bigskip

\centerline{by}

\bigskip

\medskip

\centerline{{\Large Holger~van Bargen\footnotemark[1]}}

\footnotetext[1]{Institut f\"ur Mathematik, MA 7-4, Technische Universit\"at Berlin, Stra\ss e des 17. Juni 136, D-10623 Berlin, van\underline{ }barg@math.tu-berlin.de }

\bigskip

\bigskip

{\leftskip=1truecm

\rightskip=1truecm

\baselineskip=15pt

\small

\noindent{\slshape\bfseries Abstract.} It has been shown by various authors under different assumptions that the diameter of a bounded non-trivial set $\gamma$ under the action of a stochastic flow grows linearly in time. We show that the asymptotic linear expansion speed if properly defined is deterministic i.e. we show for a $2$-dimensional isotropic Brownian flow $\Phi$ with a positive Lyapunov exponent that there exists a non-random set $\mathcal{B}$ such that we have for $\eps>0$, arbitrary connected $\gamma\subset\subset\R^2$ consisting of at least two different points and arbitrarily large times $T$ that $$(1-\eps)T\mathcal{B}\subset \cup_{0\leq t\leq T}\cup_{x\in\gamma} \Phi_{0,t}(x)\subset(1+\eps)T\mathcal{B}.$$  


\bigskip

}
\pagestyle{myheadings}\pagenumbering{arabic}
\markboth{Holger van Bargen}{Asymptotic Growth Of IBF}
\section{Introduction And Preliminaries}
Isotropic Brownian flows (IBFs) are a fairly natural class of stochastic flows and have been studied by various authors in different directions, e.g.~\cite{bh}, \cite{lj}, \cite{gd} and~\cite{vB} - just to name a few references. ~\cite{css} and \cite{ls} study the evolution of the diameter of a bounded and non-trivial set under the evolution of such a flow giving upper and lower bounds for the linear growth rate. Nevertheless these bounds turn out to be far from each other in some examples and there is little hope to match these bounds with the methods from~\cite{css} or~\cite{ls}. We will follow a different approach which first appeared in~\cite{dkk}, wherein a class of periodic stochastic flows on $\R^2$ (or stochastic flows on the torus) is considered. \cite{dkk} develop a similar limit theorem (even with a stronger assertion) using the fact that their model essentially lives on a compact manifold. Although we will sometimes follow the lines of thought of~\cite{dkk} in the first part, we will see that to get the assertion we will have to replace the methods relying on the assumption of periodicity (which means perfect dependence of particles which are far from each other) on $\R^2$ by different ones. This is done using the invariance properties with respect to time reversal of IBFs. These properties are not shared by the model of~\cite{dkk} and hence are a novelty in the present subject. The paper is divided into several sections. First we briefly review the important definitions and cite some prerequisites from the literature including a subsection on the smoothness of the density of the two-point motion. Afterwards we give the proper definition of the asymptotic linear expansion speed and state the main result, from which the fact, that the asymptotic expansion speed is constant, turns out to be a corollary. We give the proofs of the main results in the last two sections.  The first of these is dedicated to the proof of the lower bound i.e. that the expansion is sufficiently fast. Here we also identify the set $\mathcal{B}$ in terms of a stable norm (which is a concept from~\cite{dkk}). We finally finish the proof in the last section by showing that the expansion is sufficiently slow, for which it will turn out to be sufficient to show that the expansion speed is independend of the initial set. We will work in general dimension $d$ where possible. But since several important features of the proof obviously fail in higher dimensions the reader might assume that $d$ is always equal to two. 
\subsection{Stochastic Flows And Stochastic Differential Equations}
Let us first state the definition of a Brownian flow.
\begin{definition}
Let $\left(\Phi_{s,t}(x,\omega): s,t\in[0,\infty),x\in\mathbb{R}^{d},\omega\in\Omega\right)$ be a continuous $\mathbb{R}^{d}$-valued random vector field defined on a probability space $(\Omega,\mathcal{F},\mathbb{P})$. $\left(\Phi_{s,t}(x,\omega)\right)$ is called a Brownian flow of diffeomorphisms, if there is a $\mathbb{P}$-null set $N \subset\Omega$ such that we have for $\omega\in N^{C}$:
\begin{enumerate}
\item $\Phi_{s,u}(\omega)=\Phi_{t,u}(\omega)\circ\Phi_{s,t}(\omega)$ and
 $\Phi_{s,s}(\omega)=\left.\textnormal{id}\right|_{\mathbb{R}^{d}}$ for any $0\leq s,u,t<\infty$,
\item $\Phi_{s,t}(\omega): \mathbb{R}^{d}\rightarrow\mathbb{R}^{d}$ is an onto map for arbitrary $0\leq s,t<\infty$,
\item $\Phi_{s,t}(x,\omega)$ is $k$ times continuously differentiable w.r.t. $x$ for any $k$,
\item $\Phi$ is Brownian i.e. for $n\in\mathbb{N}$ and $0\leq s_{1}<s_{2}<\ldots<s_{n}<\infty$ we have that the family of random variables 
$\left(\Phi_{s_{j-1},s_{j}}:j\in\{1,\ldots,n\}\right)$ is independent.
\end{enumerate}
\end{definition}
We will write $\Phi_{t}(x,\omega),\ldots$ for $\Phi_{0,t}(x,\omega),\ldots$ and for $x,y,z,\ldots\in\mathbb{R}^{d}$ we abbriviate $x_{t}:=\Phi_{t}(x),\ldots$. Due to \cite[Theorem 4.4.1]{k} stochastic flows are generated by Kunita-type stochastic differential equations of the form
\begin{equation}dX(t)=M(dt,X(t))\label{eq:flussgleichung}\end{equation} wherein $M$ is a suitable semimartingale field. 
We will briefly describe the construction of the fields leading to isotropic Brownian flows in  the sequel. See \cite{ls} or \cite{bh} for further details. We choose a modification of $\Phi$ that satisfies the above with $N=\emptyset$.
\subsection{Covariance Tensors And Correlation Functions}
\begin{definition}
A function $b:\mathbb{R}^{d}\rightarrow\mathbb{R}^{d\times d}$ is an isotropic covariance tensor if $x\mapsto b(x)$ is $C^\infty$ and all derivatives of any order are bounded, $b(0)=E_{d}$ (the $d$-dimensional unit matrix),  $x\mapsto b(x)$ is not constant and
$b(x)=O^{*}b(Ox)O$ holds for any $x\in \mathbb{R}^{d}$ and any orthogonal matrix $O$.
\end{definition}
\textbf{Remark}: The assumptions on the differentiability of the flow as well as of the generating tensor are a bit restrictive, but we do not want to mess with smoothness problems coming especially from Malliavin calculus (we strongly conjecture that for the $2$-dimensional case a $C^6_b$-assumption should be sufficient - see~\cite{vBdiss}). 
\begin{lemma}\label{le:isoko}\label{le:korrfunk}
Let $b$ be an isotropic covariance tensor. The functions $B_L$ and $B_N$ defined by $B_{L}(r):=b_{ii}(r e_{i}), r\geq0$ and $B_{N}(r):=b_{ii}(r e_{j}), r\geq0,i\neq j$ are the longitudinal (and normal - respectively) correlation functions of $b$. 
Their definition does not depend on the specific choice of $1\leq i,j\leq d$ and we have for arbitrary $i,j\in\{1,\ldots,d\}$ and $x\in\mathbb{R}^{d}$:
\[ 
b^{ij}(x)=\left\{\begin{array}{ccc}
(B_{L}(|x|)-B_{N}(|x|))x^{i}x^{j}/|x|^{2}+B_{N}(|x|)\delta^{ij}&:&x\neq0\\
\delta^{ij}=\delta^{ij}B_{L}(0)=\delta^{ij}B_{N}(0)&:&x=0
\end{array} \right.
\]  
The right-hand derivatives of $B_{L/N}$ satisfy $\beta_{L/N}:=-B_{L/N}''(0)>0$ and we have the Taylor-expansions $B_{L/N}(r)=1-\frac{1}{2}\beta_{L/N}r^{2}+O(r^{4}):(r\rightarrow0)$ as well as the global estimate $||B_{L/N}||_{\infty}= 1$. We will use the above Taylor-expansion in the following weaker form. For any $\eps>0$ there is $r^{(\eps)}>0$ such that we have for $0<r<r^{(\eps)}$ that
\begin{equation}\label{eq:B-approximationsgleichung}
\frac{\left|1-B_{L}(r)-\frac{1}{2}\beta_{L}r^{2}\right|\vee \left|1-B_{N}(r)-\frac{1}{2}\beta_{N}r^{2}\right|}{r^{3}}<\epsilon.
\end{equation}
\end{lemma}
Proof: \cite[(2.5), (2.6), (2.13), (2.8), (2.9) and (2.18)]{bh}.\hfill$\Box$\\
\subsection{Brownian Fields And Generated Flows}
Now we can define the semimartingale field $M$ which in fact is a martingale field. 
\begin{definition}\label{def:isoko}
Let $b$ be an isotropic covariance tensor. An $\mathbb{R}^{d}$-valued random vector field $\left(M(t,x): t\geq0, x\in\mathbb{R}^{d}\right)$ - defined on a probability space $(\Omega,\mathcal{F},\mathbb{P})$ - is an isotropic Brownian field if 
$(t,x)\mapsto M(t,x)$ is a centered Gaussian process with
$\mathit{cov}(M(s,x),M(t,y))=(s\wedge t)b(x-y)$ and
$(t,x)\mapsto M(t,x)$ is continuous for almost all $\omega$.
A stochastic flow, generated via (\ref{eq:flussgleichung}) by an isotropic Brownian field is called an isotropic Brownian flow (IBF).
\end{definition}
\begin{theorem}
The $n$-point-motion $(x^{(1)}_{t},\ldots,x^{(n)}_{t}):=(\Phi_{t}\left(x^{(1)}\right),\ldots,\Phi_{t}\left(x^{(n)}\right))$ is a $\mathbb{R}^{nd}$-valued diffusion with the following properties:
\begin{enumerate}
\item For  $g\in C^{2}_{b}$ its generator $L^{(n)}g(x^{(1)},\ldots,x^{(n)})$ given  by
\[
\frac{1}{2}\sum_{l,m=1}^{n}\sum_{i,j=1}^{d}b\left(x^{(l)}-x^{(m)}\right)\frac{\partial^2g}{\partial x^{(l)}_i\partial x^{(m)}_j}\left(x^{(1)},\ldots,x^{(n)}\right).
\]
\item There is a standard Brownian motion $W$ such that $\rho^{xy}_{t}$ solves the SDE 
\begin{equation}\label{eq:zweipunktgleichung}
d\rho^{xy}_{t}=(d-1)\left(\frac{1-B_{N}(\rho^{xy}_{t})}{\rho^{xy}_{t}}\right)dt+\sqrt{2(1-B_{L}(\rho^{xy}_{t}))}dW_{t}.      
\end{equation} 
\end{enumerate}
\end{theorem}
Proof: \cite[p. 617]{lj}, \cite[p.124]{k}, \cite[p. 4]{ls} and \cite[(3.11)]{bh}.\hfill$\Box$\newline
\textbf{Remark:} \cite{lj} uses a slightly different definition. Assume $\alpha=1$ there to get things into line with the definitions above. The previous theorem shows, that $(x_{t},y_{t})$ coincides in law with the solution of the following SDE.
\begin{equation}\label{eq:zweipunktSDE}
\left(\begin{array}{c}x_{t}-x\\y_{t}-y\end{array}\right)=
\int_{0}^{t}\left(\begin{array}{cc}E_{d}&b\left(x_{s}-y_{s}\right)\\b\left(x_{s}-y_{s}\right)&E_{d}\end{array}\right)^{1/2}d{W}_{s}=:\int_{0}^{t}\bar{b}\left(x_{s}-y_{s}\right)dW_{s}
\end{equation}
Therein $W$ is a $2d$-dimensional standard Brownian motion. The following lemma states some information about the eigenvalues of $b$ and $\bar b$ respectively. 
\begin{lemma}\label{le:eigenwerte}
For $z\in\mathbb{R}^{d}$ we have:
\begin{enumerate}
\item $z$ is an eigenvector of $b(z)$ to the eigenvalue $B_{L}\left(\left|z\right|\right)$. 
\item Any vector $0\neq z^{\bot}$ perpendicular to $z$ is an eigenvector of $b(z)$ to the eigenvalue $B_{N}\left(\left|z\right|\right)$.
\item $\bar{b}$ has the eigenvalues  $\left\{1\pm B_{L}(z),1\pm B_{N}(z) \right\}$ with multiplicities $1$ and $n-1$ respectively.
\end{enumerate}
\end{lemma}
Proof: straightforward computations using Lemma~\ref{le:isoko}.\hfill$\Box$\\ Observe that the previous lemma ensures that $\bar{b}$ is elliptic apart from the diagonal $\{x=y\}$.
\subsection{Density Of The $n$-Point Motion}
As we have already seen, the one-point-motion of an IBF is a standard Brownian motion and so of course possesses a $C^\infty$-density. This section is devoted to the question if this is true for the $2$-point-motion $(x_t,y_t)$ of an $d$-dimensional IBF ($x\neq y$). The homeomorphic properties of the flow do not allow for $x_t=y_t$ to hold at any time except on a null set (remember that we decided to modify the flow in a way such that $x_t=y_t$ is impossible). One might expect the process $(x_t,y_t)$ to posesses a density on $\R^{2d}_\times:=\R^{2d}\setminus\{ z\in\R^{2d}:z_i=z_{d+i}\forall i=1,\ldots,d\}$. This is in fact true as we shall see in the following. 
\begin{theorem}\label{th:posdens}
The two-point-motion $(x_t,y_t)$ interpreted as a diffusion on $\R^{2d}_\times$ possesses a strictly positive $C^\infty$-density on $\R^{2d}_\times$.
\end{theorem}
Proof: We restrict ourselves to $t=1$ by scaling. First observe that our smoothness assumptions on $b$ allow for the use of H\"ormander's Theorem \cite{ho}. See~\cite{nu} for details and stochastic interpretations. Since we already observed that the process satisfies the SDE~\eqref{eq:zweipunktSDE} and since  Lemma~\ref{le:eigenwerte} ensures that H\"ormander's condition is satisfied we can conclude that on $\R^{2d}_\times$ a $C^\infty$-density exists. We now have to show that it is  strictly positive there. We want to apply the results of~\cite{l1}, so we have to consider the following control problem.
\be
dz_t(h)=\bar{b}(z_t(h))h_t dt
\ee 
Therein $h$ is a square-integrable, $\R^{2d}$-valued control function (in fact chosen to be continuously differentiable). $z_t$ is a $2d$-dimensional process to be thought of as a deterministic version of the two-point-motion. Fix $(x,y)\in \R^{2d}_\times$. In  order to show that $(x_1,y_1)$ has positive transition density for any $x^{(1)},y^{(1)}\in \R^{2d}_\times$ it is enough to establish the following Bismut Condition (see~\cite{bi})
\begin{condition}
For any $(x,y)=z\in \R^{2d}_\times,(x^{(1)},y^{(1)})\in \R^{2d}_\times$ there is an $h\in L^2$ such that 
\be z_1(h)=(x^{(1)},y^{(1)}) \label{eq:ankommen}\ee and such that $h\mapsto(z_1(h))$ is a submersion in $h$. (we identify $\R^{2d}$ and $\R^{d}\times\R^d$ in the obvious way). 
\end{condition}
Proof.: Step 1: Let us assume first that $\overline{x,x^{(1)}}$ and $\overline{y,y^{(1)}}$ are disjoint and that each of them consists at least of two points. ($\overline{x,y}$ denoting the convex hull of $x$ and $y$.) We construct a control satisfying~\eqref{eq:ankommen} such that the stream lines of $z_t$ are exactly $\overline{x,x^{(1)}}\cup\overline{y,y^{(1)}}$ . This ensures that $\bar{b}(z_t(h))$ is regular and its determinant is bounded away from zero for all $t$. The simplest way to obtain the desired streamlines is to ensure $\bar{b}(z_t(h))h_t\equiv \left(z^{(1)}-\begin{pmatrix}x\\y\end{pmatrix}\right)$. We may hope to achieve this by setting $h_0:=\bar{b}\left(\begin{pmatrix}x\\y\end{pmatrix}\right)^{-1}\left(z^{(1)}-\begin{pmatrix}x\\y\end{pmatrix}\right)$ as well as $0=\frac{d}{dt}\left[\bar{b}(z_t(h)) h_t\right] $ which is the same as
\be \label{eq:controlode} 
\bar{b}(z_t(h))\frac{dh_t}{dt}=-\left[\left( \ip{\frac{dz_t(h)}{dt}}{\nabla} \bar{b}\right)(z_t(h))\right]h_t
\ee
So we see that we can choose $h_t$ to be the projection on the first $2d$ coordinates of the solution to the following $4d$-dimensional initial value problem.
 \be
\left\{\begin{array}{rl}
   \frac{d}{dt}\begin{pmatrix}h_t\\z_t(h) \end{pmatrix}=&\begin{pmatrix}-\bar{b}^{-1}(z_t(h))\left[\left(\ip{\frac{dz_t(h)}{dt}}{\nabla} \bar{b}\right)(z_t(h))\right]h_t\\ \bar{b}(z_t(h))h_t \end{pmatrix},\\
h_0=&\bar{b}\left(\begin{pmatrix}x\\y\end{pmatrix}\right)^{-1}\left(z^{(1)}-\begin{pmatrix}x\\y\end{pmatrix}\right),\hspace{1cm}z_0(h)=\begin{pmatrix}x\\y \end{pmatrix}
\end{array}\right.
\ee
Existence and uniqueness of a solution to this initial value problem can be obtained from the standard theorems because we ensured that the determinant of $\bar{b}$ is bounded away from zero and hence that the right-hand-side of \eqref{eq:controlode} is continuously differentiable. 
\\Step 2: For a general positions of $x$, $x^{(1)}$, $y$ and $y^{(1)}$ observe that we can divide the action into two parts i.e. timesteps of length 0.5 and choose the streamlines of $x$ and $y$ to be piecewise linear and disjoint.\\ 
Step 3: Finally we have to note that by Theorems 1.1 (smoothness) and 1.10 (surjectivity) of~\cite{bi} we have a submersion in $h$. \hfill$\Box$\\
\subsection{Lyapunov-Exponents}
As proved in \cite[(7.2) and (7.3)]{bh} IBFs have Lyapunov exponents which satisfy $\mu_{i}:=\frac{1}{2}\left[(d-i)\beta_{N}-i\beta_{L}\right]$. 
The top Lyapunov-exponent $\mu_1$ i.e. its sign crucially affects the asymptotic behaviour of the flow, as shown in~\cite{css}. We make the standing assumption $\mu_{1}=\frac{1}{2}[(d-1)\beta_{N}-\beta_{L}]>0$. If this is not fulfilled then~\cite{ss} shows that our main result cannot be expected to be true because the flow contracts a closed ball of positive diameter to a point with positive probability. 
\subsection{Support Theorem For Isotropic Brownian Flows}
As to every Gaussian measure one can associate a Hilbert space to an isotropic Brownian field - the so-called reproducing kernel Hilbert space $\mathcal{H}$ of $M$. For details on this we refer the reader to~\cite{b} in the general context of Gaussian measures and to~\cite{gd} the special case considered here. We will only need the fact from~\cite{gd} that for $x\in\Rd$ and an arbitrary signed measure $\mu$ on the Borel sets of $\mathbb{R}^{d}$ the vectorfield $\int b^{i,j}(x-y)d\mu(y,j)$ belongs to $\mathcal{H}$. 
\begin{theorem}\label{sa:ibfsupp}
Let $M$ be an isotropic Brownian field. Due to~\cite[Theorem 3.5.1]{b}  this can be written as $M(t,x)=\sum_{i=1}^{\infty}V_{i}(x)W_{t}^{i}$ wherein $(V_i)_{i\in\N}$ is a complete orthonormal system in $\mathcal{H}$. Assume that $V_{1}$ is four times continuously differentiable and that all derivatives up to order four are bounded. Then for $\mathbb{K}\subset\subset\mathbb{R}^{d}$, $T>0$ and $\delta>0$ there are positive numbers $\epsilon$ and $C_1$, such that: 
\[
\prob{\sup_{0\leq t\leq T}\sup_{x\in\mathbb{K}}\left\|x_{\frac{t}{C_1}}-\psi_{t}(x) \right\|<\delta}>\epsilon
\]
Therein $\psi=\psi_{t}(x)$ is the solution of the following deterministic control problem:
\[
\left\{
\begin{array}{ccc}
 \partial_{t}\psi_{t}(x)&=&V_{1}(\psi_{t}(x))\\
\psi_{0}(x)&=&x
\end{array}
 \right.
\] 
\end{theorem}
Proof: This is Theorem 6.2.3 of~\cite{gd}.
\hfill$\Box$\newline
For the convenience of the reader we include the following definition.
\begin{definition} An $\mathcal{H}$-simple control $V$ is a mapping from $[0,T]$ to $\mathcal{H}$, which is piecewise constant.
\end{definition}
\subsection{Time Reverse And Markov-Properties}
\begin{lemma}\label{le:reverse}
For arbitrary $T>0$ we have
\begin{equation}
\mathcal{L}\left[\left(\Phi_{s,t}(.):0\leq s\leq t\leq T \right)\right]=\mathcal{L}\left[\left(\Phi_{T-s,T-t}(.):0\leq s\leq t\leq T \right)\right].
\end{equation}
\end{lemma}
Proof: Due to \cite[Theorem 4.2.10]{k} the backward flow is driven by the same infinitesimal generator as the forward flow (see Lemma~\ref{le:korrfunk} for details). Therefore the law of the forward flow and the law of the backward flow coincide.\hfill$\Box$\\
\begin{lemma}
Let $\mathcal{F}_{s,t}$ be the $\sigma$-field generated by $\left\{\Phi_{r,u}:s\leq r\leq u \leq t \right\}$ .
\begin{enumerate}
\item For an $\left(\mathcal{F}_{s,t}:t\in[s,\infty) \right)$-stopping-time $\tau$ we have:
\begin{equation}
\mathcal{L}\left[\left(\Phi_{\tau,r}\left(\Phi_{s,\tau}(.) \right):r\geq\tau\right) \right|\left.\mathcal{F}_{s,\tau}\right]=\mathcal{L}_{s,\tau}\left[\Phi_{s,s+r-\tau}(.):r\geq\tau\right]
\end{equation}
\item For any $\left\{\mathcal{F}_{s,t}:s\in(-\infty,t]\right\}$-stopping-time $\tau$ we have:
\begin{equation}
\mathcal{L}\left[\Phi_{\tau,t}\left(\Phi_{r,\tau}(.) \right):r\leq\tau \right|\left.\mathcal{F}_{\tau,t}\right]=\mathcal{L}_{\tau,t}\left[\Phi_{t+r-\tau,t}(.):r\leq\tau\right]
\end{equation}
\end{enumerate}
\end{lemma}
Proof: 2. is a consequence of 1. and Lemma~\ref{le:reverse}. 1. is~\cite[Theorem~4.2.1]{k}.\hfill $\Box$\\
\subsection{Chasing Ball Property, LDP For Discrete Supermartingales}
The first of the following lemmas states that the distance of a non-trivial set under the action of the flow tends to approach another moving particle (arbitrary non anticipating movement) provided that the other particle does not move too fast. Therein we mean  
\begin{definition}A subset of $\mathbb{R}^{d}$ is called non-trivial, if it is bounded, connected and consists of at least two different points.\end{definition}
Note that for IBFs the estimates of the local characteristics and the ellipticity bounds of~\cite{ss} hold. Therefore we may use the following lemma.
For $t\geq0$ denote by $\mathcal{F}_{t}:=\mathcal{F}_{0,t}$ the sigma-field, generated by the flow up to time $t$. 
\begin{lemma}\label{le:lokdrift}
Let $\Phi$ be an IBF with generator $M$. Then there are functions $G':[0,\infty)\times[0,\infty)\times[0,\infty)\rightarrow[0,\infty)$ and $G'':[0,\infty)\times[0,\infty)\rightarrow[0,\infty)$, such that there is $r_0>0$ depending only on $b$ such that we have the following.
\begin{enumerate}
\item For all $s\in [0,\infty)$ the function $G'(\cdot,s,\cdot)$ is continuous, non-increasing with\\ $\lim_{K\rightarrow\infty}\lim_{r\rightarrow\infty}G'(K,s,r)=0$.
\item For all $s\in [0,\infty)$  $G''(s,\cdot)$ is continuous and $r\in(0,r_0)\Rightarrow G''(s,r)>0$.
\item Let $s>0$ and $r<r_0$. Let $\tau$ be a finite stopping time for the flow and $x,y,z$ $\mathcal{F}_{\tau}$-measurable random points in $\mathbb{R}^{d}$ with $\left\|x-y\right\|=r$. 
Define $r_{1}:=\left\|x-z\right\|\wedge\left\|y-z\right\|$,  $r_{2}:=\left\|\Phi_{\tau,\tau+s}(x)-z\right\|\wedge\left\|\Phi_{\tau,\tau+s}(y)-z\right\|$. 
Then we have \[\ce{r_{2}\vee (r_1-K)}{\mathcal{F}_{\tau}}\leq r_{1}+G'(K,s,r_{1})-G''(s,r).\] 

\end{enumerate}
\end{lemma}
Proof:~\cite[Lemma 2.5]{ss}. Observe that $K$ does not appear in the original result in~\cite{ss} but can be obtained by adding it in the proof of (15) on pages 2055 and 2056 of~\cite{ss} to obtain instead of (15) the estimate $\expec{||x_{\tau+s}||\vee (x^1+K)}-x^1\leq g(K,x^1)$ with $\lim_{x\to\infty}g(K,x)=g_K$ and $\lim_{K\to\infty}g_K=0$ and by proceeding as in~\cite{ss} afterwards.\hfill$\Box$\\
The next lemma is an elementary large deviation principle and we recall it for the convenience of the reader.
 \begin{lemma}\label{le:mmb}
Let $\{\xi_{j}:j\in\mathbb{N}\}$ be a sequence of real-valued random variables with
\begin{enumerate}
\item $\ce{\xi_{j+1}}{\xi_{1},\ldots,\xi_{j}}\leq0$,
\item $\forall m\in\mathbb{N}:\exists K_{m}\in\mathbb{R}:\forall j\in\mathbb{N}:\mathbb{E}\left[|\xi_{j}|^{m}\right]\leq K_{m}$.
\end{enumerate}
Then we have that for $\eps>0$ there exist constants $\kappa^{(1)}_m$ depending on $\epsilon$ and $(K_{n})_{n\in\N}$ such that we have for $n\in\N$ that
$\mathbb{P}\left[\sum_{j=1}^{n}\xi_{j}\geq\epsilon n\right]\leq\kappa^{(1)}_m n^{-m}$.
\end{lemma}
Proof: \cite[Lemma 2]{dkk}.\hfill$\Box$\\
\subsection{Sub-Gaussian Tails And Sublinear Growth}
\begin{lemma}\label{le:subgauss}
There is a positive constant $C_2$, such that $\mathbb{P}$-a.s. for any bounded subset $\gamma$ of $\mathbb{R}^{d}$ we have
$\limsup_{T\rightarrow\infty}(\sup_{t\in[0,T]}\sup_{x\in\gamma}\frac{1}{T}\left\|\Phi_{t}(x)\right\|)\leq C_2.$
\end{lemma}
Proof: \cite[Theorem 2.1]{ls2}.\hfill$\Box$\\
\section{Statement Of The Main Results}
\begin{theorem}\label{th:main}
For any bounded, connected $\gamma\subset\R^2$ consisting of at least two different points we let $\gamma_t:=\Phi_t(\gamma) $ and $\mathcal{W}_{t}(\gamma):=\bigcup_{0\leq s\leq t}\gamma_{s}$. Then there exists a deterministic set $\mathcal{B}$ such that we get for any $\eps>0$:
\begin{enumerate}
\item There is $\mathbb{P}$-a.s. $0<T(\gamma,\epsilon)<\infty$, such that for any $t>T(\gamma,\epsilon)$  the following holds.
\begin{equation}(1-\epsilon)t\mathcal{B}\subset\mathcal{W}_{t}(\gamma).\end{equation}
\item There is $\mathbb{P}$-a.s. a sequence $\left(t_{k}:k\in\mathbb{N}\right)\subset\R_+$ with $t_{k}\nearrow\infty$ that fulfills \be \mathcal{W}_{t_{k}}(\gamma)\subset(1+\epsilon)t_{k}\mathcal{B}.\ee
\end{enumerate}
We also have 
\be \lim_{T\to\infty}\prob{(1-\epsilon)T\mathcal{B}\subset\mathcal{W}_T(\gamma)\subset(1+\epsilon)T\mathcal{B}}=1.
\ee
\end{theorem} 
Proof: The proof will be given in the sections~\ref{sec:lb} and~\ref{sec:ub}.
\begin{corollary}
If we define for $\gamma$ as above the asymptotic linear expansion speed to be 
$$\liminf_{T\to\infty}\frac{1}{T}\diam(\mathcal{W}_T(\gamma))$$ then it is independend of $\gamma$ and a.s. constant.
\end{corollary}
Proof: This follows directly from Theorem~\ref{th:main}.\hfill$\Box$\\
\section{The Lower Bound}
\label{sec:lb}
\subsection{Hitting Time Of Far Away Balls}
Assume that the original set $\gamma\subset\mathbb{R}^{d}$ is connected, compact and that it consists of at least two different points (the assumption of compactness is made for simplicity and could be omitted). Denote by $\gamma_{t}:=\Phi_{t}(\gamma)$ the set $\gamma$ at time $t$ and by $d_{t}:=\diam(\gamma_{t})$ its diameter. Further denote for any $R>0$ by $\tau^{R}(\gamma,P):=\inf\left\{t>0:\dist(\gamma_{t},P)\leq R, d_{t}\geq 1\right\}$ the time it takes for $\gamma$ to reach an $R$-neighbourhood of $P\in\Rd$. In fact it will turn out that $\liminf_{t\rightarrow\infty}d_{t}\geq 1$ $a.s.$. We call a subset of $\mathbb{R}^{d}$ large if it is bounded and has diameter at least $1$.  
Due to the results of \cite{ss} and \cite{css} we may assume that $\gamma$ is large (the following will proof that $\gamma$ will become large a.s. anyways).
\begin{theorem}\label{th:cor6}
Let $P\in\Rd$, assume that $\gamma\subset\Rd$ is large and define $\bar{r}:=1\vee\mathit{dist}(P,\gamma)$. There is a constant $R>0$ such that for any
$m\in\mathbb{N}$ there is $\kappa^{(2)}_{m}>0$ (not depending on $\gamma$, $P$ or $\bar{r}$) such that for $\beta>1$ we have 
$\mathbb{P}\left[\tau^{R}(\gamma,P)>\kappa^{(2)}_m\beta \bar{r}\right]\leq \kappa^{(2)}_{m}\beta^{-m}\bar{r}^{-m}$.
\end{theorem}
The proof consists of several steps:
\begin{enumerate}
\item Construction of a strictly increasing $C^{2}$-function  $f:(0,\infty)\rightarrow\mathbb{R}$  with $\lim_{r\rightarrow\infty}f(r)=\infty$ such that $f(\rho^{xy}_{t})$ is a submartingale for any $x, y\in\Rd $. The drift of this submartingale has to be bounded away from zero for small $\rho_{t}^{xy}$.
\item Getting estimates on the growth of $d_{t}$ on average
\item Getting estimates of the probability of finding $\gamma_{t}$ not being large after a long time
\item Establishing a negative upper bound for the \glqq drift\grqq of $r_{t}:=\dist(\gamma_{t},P)$ outside $K_{R}(P):=\left\{x\in\Rd:|x-P|\leq R \right\}$ 
\item Glue all the above together to prove Theorem~\ref{th:cor6}
\end{enumerate}
\subsubsection{Construction Of $f$}
The first ingredient needed to construct $f$ is the following lemma. 
\begin{lemma}
For any $0<c_8<c_9<\infty$, $\delta>0$ and $-\infty<c_{10}<0<c_{11}<\infty$  there is a decreasing $C^{2}$-function $h:[c_8,c_9]\rightarrow [h(c_9),0]$ with \label{le:ableitungskorrektur}
\begin{enumerate}
\item $h'(c_8)=h'(c_9)=0$,
\item $h''(c_8)=c_{10}$, $h''(c_9)=c_{11}$ and $h''$ is increasing,
\item $h(c_8)=0$,
\item $\sup_{c_8\leq r\leq c_9}\{|h'(r)|\}\leq\delta$. 
\end{enumerate}
\end{lemma}
Proof: For $0<\epsilon<0,5(c_{11}\wedge-c_{10})$ define $h''_{\epsilon}:[c_8,c_9]\rightarrow[c_{10},c_{11}]$ via\\ $h''_{\epsilon}(r):=\ind{\left[c_8,c_{8,\epsilon}\right]}(r)(r-c_{8,\epsilon})\frac{c_{10}^{2}}{2\epsilon(c_9-c_8)}+\ind{\left[c_{9,\epsilon},c_9\right]}(r)(r-c_{9,\epsilon})\frac{c_{11}^{2}}{2\epsilon(c_9-c_8)}$ with\\ $c_{8,\epsilon}:=c_8-\frac{2\epsilon(c_9-c_8)}{c_{10}}$ and $c_{9,\epsilon}:=c_9-\frac{2\epsilon(c_9-c_8)}{c_{11}}$.
Letting also
$ h'_{\epsilon}(r):=\int_{c_8}^{r}h''_{\epsilon}(s)ds=\ind{[c_{9,\epsilon},c_9]}(r)(r-c_{9,\epsilon})^{2}\frac{c_{11}^2}{4\epsilon(c_9-c_8)}-\ind{]c_{8,\epsilon},c_9]}(r)\epsilon(c_9-c_8)$\\
$\hspace{5mm}+\ind{[c_8,c_{8,\epsilon}]}(r)\left[(r-c_{8,\epsilon})^{2}-(c_8-c_{8,\epsilon})^{2}\right]\frac{c_{10}^2}{4\epsilon(c_9-c_8)}$  ensures 1.. We then also have  $h'_{\epsilon}\leq0$. 4. follows from choosing $\epsilon\leq\frac{\delta}{c_9-c_8}$.
Setting $h(r):=\int_{c_8}^{r}h'_{\epsilon}(s)ds$ for such an $\epsilon$ makes $h$ decreasing, ensures 3. and finishes the proof of Lemma~\ref{le:ableitungskorrektur}.
\hfill$\Box$\newline
\begin{lemma}\label{le:exf}
There is a strictly increasing $C^2$-function $f$ with the following properties.
\begin{enumerate}
\item $\lim_{r\rightarrow\infty}f(r)=\infty$ and $f(\rho^{xy}_{t})$ is a submartingale for any $x, y\in\Rd $
\item $f(1)=0$
\item Writing with It\^{o}s formula, (\ref{eq:zweipunktgleichung}) and Fubinis theorem
\bee  
&&\mathbb{E}\left[f(\rho^{xy}_{t+s})-f(\rho^{xy}_{s})\right]\nonumber\\ &=&\int_{s}^{t+s}\mathbb{E}\left[f'(\rho^{xy}_{r})\frac{1-B_{N}(\rho^{xy}_{r})}{\rho^{xy}_{r}}(d-1)+f''(\rho^{xy}_{r})\left(1-B_{L}(\rho^{xy}_{r})\right)\right]dr\nonumber\\
&+&\mathbb{E}\left[\int_{s}^{t+s}f'(\rho^{xy}_{r})\sqrt{2(1-B_{L}(\rho^{xy}_{r}))}dW_{r}\right]\nonumber\\
&=:&\int_{s}^{t+s}\mathbb{E}\left[g(\rho^{xy}_{r})\right]dr
+\mathbb{E}\left[\int_{s}^{t+s}\tilde{g}(\rho^{xy}_{r})dW_{r}\right]
\eee   we get that $\tilde{g}$ is bounded and that $g-\frac{1}{8}(\beta_{N}(d-1)-\beta_{L})\geq0$.
\item There are $C_8>0$ and $C_9>0$ such that 
\begin{equation}\ce{(f(d_{t+1})-f(d_{t}))\wedge C_9}{\mathcal{F}_{t}}\geq C_7.\label{eq:driftschranke}\end{equation}
($\mathcal{F}_{t}$ denotes the $\sigma$-field generated by the flow up to time $t$.)
\end{enumerate}
\end{lemma}
Proof: We choose the following ansatz for $f$ which uses a local linearization of (\ref{eq:zweipunktgleichung}) near the origin. 
\begin{equation}
f(r)-c_{1}:=\left\{\begin{array}{ll}
\log r+c_{2}&:0<r<c_8\\
c_{3}\sqrt{r}+h(r)&:c_8\leq r\leq c_9\\
c_{4}r+c_{5}&:c_9<r
\end{array}\right.
\end{equation}Put
\begin{align} \epsilon:=&1\wedge\frac{1}{8}\frac{\beta_{N}(d-1)-\beta_{L}}{d-1}\wedge\frac{1}{24}(\beta_{N}(d-1)-\beta_{L})\left(\frac{\beta_{L}}{\beta_{N}(d-1)}\right)^{\frac{1}{3}}\nonumber\\
&\wedge\frac{1}{24}(\beta_{N}(d-1)-\beta_{L}) \left(\frac{\beta_{N}(d-1)}{\beta_{L}}\right)^{-\frac{4}{3}} \end{align}   and choose $r_{\epsilon}$ according to (\ref{eq:B-approximationsgleichung}).
Further set
\begin{align}
c_9:=&r_{\epsilon}\wedge1,\hspace{1cm}c_8:=c_9\left(\frac{\beta_{N}(d-1)}{\beta_{L}}\right)^{-\frac{2}{3}},\nonumber\\
\delta:=&\frac{1}{24}(\beta_{N}(d-1)-\beta_{L})\left\|\frac{1-B_{N}(.)}{c_8}\right\|_{\infty}^{-1}(d-1)^{-1},\nonumber\\
c_{11}:=&\frac{1}{2c_9\sqrt{c_8 c_9}},\hspace{0.5cm}\textnormal{ and }\hspace{0.5cm}c_{10}:=-\frac{1}{2c_8^{2}}.\nonumber
\end{align}
Choosing $h$ according to Lemma~\ref{le:ableitungskorrektur} and
\begin{align}
c_{3}:=&\frac{2}{\sqrt{c_8}},\hspace{1cm}c_{2}:=2-\log(c_8),\hspace{1cm}c_{4}:=\frac{c_{3}}{2\sqrt{c_9}}+h'(c_9)=\frac{1}{\sqrt{c_8 c_9}},\nonumber\\
c_{5}:=&\sqrt{\frac{c_9}{c_8}}+h(c_9),\hspace{1cm}c_{1}:=-c_{4}-c_{5}\nonumber
\end{align}
ensures the $C^2$-property of $f$. Its submartingale property will follow if we can show that $g$ is strictly positive and that $\tilde{g}$ is bounded. To check this let us give $f$ and its derivatives in terms of $c_8$, $c_9$ and $\delta$.
\begin{align}
f(r)+\sqrt{\frac{c_9}{c_8}}+h(c_9)+\frac{1}{\sqrt{c_8c_9}}&=&\left\{\begin{array}{ll}
\log r+2-\log(c_8)&:0<r<c_8\\
\frac{2}{\sqrt{c_8}}\sqrt{r}+h(r)&:c_8\leq r\leq c_9\\
\frac{r}{\sqrt{c_8c_9}}+\sqrt{\frac{c_9}{c_8}}+h(c_9)&:c_9<r
\end{array}\right.\nonumber
,\end{align}
\[ \left(f'(r),f''(r)\right)=\left\{\begin{array}{ll}
\left(\frac{1}{r},-\frac{1}{r^{2}}\right)&:0<r<c_8\\
\left(\frac{1}{\sqrt{c_8r}}+h'(r),-\frac{1}{2r\sqrt{c_8r}}+h''(r)\right)&:c_8\leq r\leq c_9\\
\left(\frac{1}{\sqrt{c_8c_9}},0\right)&:c_9<r
\end{array}\right.
...\]
The computation for the boundedness $\tilde{g}$ is rather simple and yields
$|\tilde{g}(r)|\leq \sqrt{\beta_{L}}+\sqrt{2}+\left|f'(c_8)\right|\left\|\sqrt{2(1-B_{L}(.))}\right\|_{\infty}<\infty$. We leave the details to the reader. Now we can turn to the estimation of $g(r)$:
For $r\geq c_9$ we obviously have $g(r)>0$ since $f'(r)>0$, $f''(r)=0$ and $B_{N}(r)<1$. 
The case $c_8\leq r\leq c_9$ needs a little more attention.
\begin{align} 
g(r)=&\frac{\beta_{N}r}{2\sqrt{c_8r}}(d-1)-\frac{\beta_{L}r^{2}}{2}\left(\frac{1}{2r\sqrt{c_8r}}-h''(r)\right)+\frac{d-1}{\sqrt{c_8r}}\left(\frac{1-B_{N}(r)-\frac{1}{2}\beta_{N}r^{2}}{r}\right)\nonumber\\
&+\left(h''(r)-\frac{1}{2r\sqrt{c_8r}}\right)\left(1-B_{L}(r)-\frac{1}{2}\beta_{L}r^{2}\right)+h'(r)\frac{1-B_{N}(r)}{r}(d-1)\nonumber\\
\geq&\frac{\beta_{N}r}{2\sqrt{c_8r}}(d-1)-\frac{\beta_{L}r^{2}}{2}\left(\frac{1}{2r\sqrt{c_8r}}+\frac{1}{2c_8^2}\right)-\frac{d-1}{\sqrt{c_8r}}\left|\frac{1-B_{N}(r)-\frac{1}{2}\beta_{N}r^{2}}{r}\right|\nonumber\\
&-\left|h''(r)-\frac{1}{2r\sqrt{c_8r}}\right|\left|1-B_{L}(r)-\frac{1}{2}\beta_{L}r^{2}\right|-\left|h'(r)\frac{1-B_{N}(r)}{r}(d-1)\right|\nonumber\\
=:&\mathit{I}-\mathit{II}-\mathit{III}-\mathit{IV}.
\end{align}
For we chose $c_8:=c_9\left(\frac{\beta_{N}(d-1)}{\beta_{L}}\right)^{-\frac{2}{3}}$ and $c_9:=r_{\epsilon}\wedge1$ we get with $\delta$ as defined above
\begin{align}
I=&\frac{\sqrt{r}}{2\sqrt{c_8}}\left[\frac{1}{2}(\beta_{N}(d-1)-\beta_{L})+\frac{1}{2}\left(\beta_{N}(d-1)-r^{\frac{3}{2}}c_8^{-\frac{3}{2}}\beta_{L}\right)\right]\nonumber\\
\geq&\frac{1}{4}(\beta_{N}(d-1)-\beta_{L})+\frac{1}{4}\sqrt{\frac{r}{c_8}}\left[\beta_{N}(d-1)-c_9^{\frac{3}{2}}c_8^{-\frac{3}{2}}\beta_{L}\right]\nonumber\\
=&\frac{1}{4}(\beta_{N}(d-1)-\beta_{L}),\\
\mathit{II}\leq&\frac{r\sqrt{rc_9}}{\sqrt{c_8r}}\left|\frac{1-B_{N}(r)-\frac{1}{2}\beta_{N}r^{2}}{r^{3}}(d-1)\right|
\leq\frac{1}{24}(\beta_{N}(d-1)-\beta_{L}),
\end{align}\begin{align}
\mathit{III}=&\left|h''(r)-\frac{1}{2r\sqrt{c_8r}}\right|\left|1-B_{L}(r)-\frac{1}{2}\beta_{L}r^{2}\right|\nonumber\\
\leq&\left|\frac{1}{2c_8^2}+\frac{1}{2r\sqrt{c_8r}}\right|r^{3}\left|\frac{1-B_{L}(r)-\frac{1}{2}\beta_{L}r^{2}}{r^{3}}\right|\nonumber\\
\leq&\frac{c_9^2}{c_8^2}r \left|\frac{1-B_{L}(r)-\frac{1}{2}\beta_{L}r^{2}}{r^{3}}\right|
\leq\frac{1}{24}(\beta_{N}(d-1)-\beta_{L}),\\
\mathit{IV}=&\left|h'(r)\frac{1-B_{N}(r)}{r}(d-1)\right|\leq\frac{1}{24}(\beta_{N}(d-1)-\beta_{L})
\end{align}
and so finally $g(r)\geq \frac{1}{8}(\beta_{N}(d-1)-\beta_{L})>0$. 
The remaining case $r\leq c_8$ similar to the above but simpler. 
It remains to show \eqref{eq:driftschranke}. This is done in the following subsubsection.
\subsubsection{Growth Of $f(d_{t})$ On Average}\label{suse:mitwachs}
\begin{figure}\centering\resizebox{10cm}{!}{\includegraphics{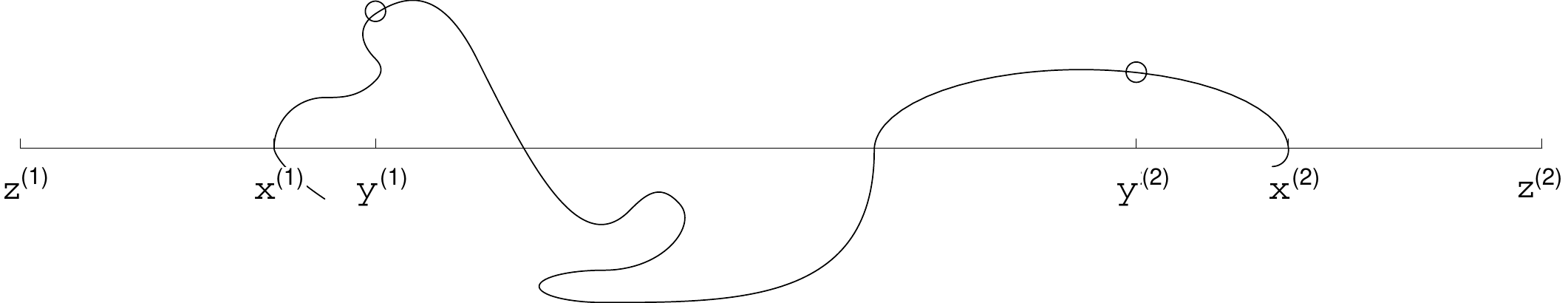}}\caption{ growth of $d_t$ on average\label{abb:mitwachs}}\end{figure}

There are two cases. If $d_{t}< \frac{r^{(\epsilon)}}{2}$ it is sufficient to consider the two-point motion. Due to the Markov-property of the submartingale $f(\rho^{xy}_{t})$ we get choosing $x_{t}$ and $y_{t}$ with $|x_{t}-y_{t}|=d_{t}$ and some constant $C_9>0$ (to be specified later)
\begin{align}
&\ce{\left((f(d_{t+1})-f(d_{t}))\wedge C_9\right)\indb{d_{t}<\frac{r^{(\epsilon)}}{2}}}{\rule{0mm}{4mm}\mathcal{F}_{t}}\nonumber\\
\geq&\ce{(f(\rho^{xy}_{t+1})-f(\rho^{xy}_{t}))\wedge C_9}{\mathcal{F}_{t}}\indb{f(d_{t})<f\left(\frac{r^{(\epsilon)}}{2}\right)}\nonumber\\
=&\mathbb{E}_{f(\rho^{xy}_{t})}\left[(f(\rho^{xy}_{1})-f(\rho^{xy}_{0}))\wedge C_9 \right]\indb{f(d_{t})<f\left(\frac{r^{(\epsilon)}}{2}\right)}\nonumber\\
\geq&\left[\left(f\left(r^{(\epsilon)}\right)-f\left(\frac{r^{(\epsilon)}}{2}\right)\right)\wedge C_9\right]\mathbb{P}_{f(\rho^{xy}_{t})}\left[\sup_{0\leq s\leq1}f(\rho^{xy}_{s})\geq f(r^{(\epsilon)})\right]\indb{f(d_{t})<f\left(\frac{r^{(\epsilon)}}{2}\right)}\nonumber\\
&+\left(\frac{1}{8}(\beta_{N}(d-1)-\beta_{L})\wedge C_9\right)\mathbb{P}_{f(\rho^{xy}_{t})}\left[\sup_{0\leq s\leq1}f(\rho^{xy}_{s})< f(r^{(\epsilon)})\right]\indb{f(d_{t})<f\left(\frac{r^{(\epsilon)}}{2}\right)}\nonumber\\
\geq&\left(\left[f\left(r^{(\epsilon)}\right)-f\left(\frac{r^{(\epsilon)}}{2}\right)\right]\wedge\frac{1}{8}(\beta_{N}(d-1)-\beta_{L})\wedge C_9\right)\indb{d_{t}<\frac{r^{(\epsilon)}}{2}}.\nonumber
\end{align} 
If $d_{t}\geq \frac{r^{(\epsilon)}}{2}$ first consider the growth of $d_{t}$. We may assume  $G''(1,\frac{r^{(\epsilon)}}{10})=:\frac{c_6}{c_4}>0$ (otherwise we decrease $r^{(\epsilon)}$, see Lemma~\ref{le:lokdrift}). There are $\hat{r}$ and $C_9>0$ such that for any $r\geq\hat{r}$ we have $G'(\frac{C_9}{2c_4},1,r)<\frac{c_6}{2c_4}$. Choose $x^{(1)}_{t},x^{(2)}_{t},y^{(1)}_{t},y^{(2)}_{t}\in\gamma_{t}$ with $|x^{(1)}_{t}-x^{(2)}_{t}|=d_{t}, |x^{(i)}_{t}-y^{(i)}_{t}|=\frac{r^{(\epsilon)}}{10}: i=1,2$ and define 
\begin{align}
z^{(1)}:=&x_{t}^{(1)}+\frac{x^{(1)}_{t}-x^{(2)}_{t}}{|x^{(1)}_{t}-x^{(2)}_{t}|}\hat{r},  z^{(2)}:=x_{t}^{(2)}+\frac{x^{(2)}_{t}-x^{(1)}_{t}}{|x^{(2)}_{t}-x^{(1)}_{t}|}\hat{r},\nonumber\\
r^{(i)}_{1}:=&|x_{t}^{(i)}-z^{(i)}|\wedge|y_{t}^{(i)}-z^{(i)}|=\hat{r}: i=1,2\textnormal{ and}\nonumber\\
r^{(i)}_{2}:=&|x_{t+1}^{(i)}-z^{(i)}|\wedge|y_{t+1}^{(i)}-z^{(i)}|:i=1,2\nonumber
\end{align}
(see Fig.~\ref{abb:mitwachs} for the geometry at time $t$).  
Lemma~\ref{le:lokdrift} provides for $i=1,2$ that we have 
$\cel{(c_4(r^{(i)}_{1}-r^{(i)}_{2}))\wedge \frac{C_9}{2}}{\mathcal{F}_{t}}\geq G''(1,\frac{r^{(\epsilon)}}{10})-G'\left(\frac{C_9}{2c_4},1,\hat{r}\right)\geq \frac{c_{6}}{2}>0$ and therefore 
$|z^{(1)}-z^{(2)}|=r^{(1)}_{1}+d_{t}+r^{(2)}_{1}$ and $|z^{(1)}-z^{(2)}|\leq r^{(1)}_{2}+d_{t+1}+r^{(2)}_{2}\Rightarrow$\\ 
$(c_4(d_{t+1}-d_{t}))\wedge C_9\geq (c_4(r^{(1)}_{1}-r^{(1)}_{2}))\wedge \frac{C_9}{2}+(c_4(r^{(2)}_{1}-r^{(2)}_{2}))\wedge \frac{C_9}{2}$ impliing\\
\begin{align}
&\ce{(c_4(d_{t+1}-d_{t}))\wedge C_9}{\mathcal{F}_{t}}\nonumber\\
\geq&\cel{(c_4(r^{(1)}_{1}-r^{(1)}_{2}))\wedge \frac{C_9}{2}}{\mathcal{F}_{t}}+\cel{(c_4(r^{(2)}_{1}-r^{(2)}_{2}))\wedge \frac{C_9}{2}}{\mathcal{F}_{t}}\geq c_{6}\nonumber.
\end{align}
Now we turn this into an estimate for $f(d_{t})$. Abbreviate $\rho_{t}:=\rho^{x^{(1)}x^{(2)}}_{t}$ and consider for $K>0$
\begin{align}\label{eq:c7def}
&\cel{\left((f(d_{t+1})-f(d_{t}))\wedge C_9\right)\indb{f\left(\frac{r^{(\epsilon)}}{2}\right)\leq f(d_{t})\leq K}}{\rule{0mm}{4mm}\mathcal{F}_{t}}\nonumber\\
\geq& \cel{\left((f\left(\rho_{t+1}\right)-f\left(\rho_{t}\right))\wedge C_9\right)\indb{f\left(\frac{r^{(\epsilon)}}{2}\right)\leq f(d_{t})\leq K}} {\rule{0mm}{4mm}\mathcal{F}_{t}}\\
\geq&\inf_{f\left(\frac{r^{(\epsilon)}}{2}\right)\leq f\left(\rho_{t}\right)\leq K}\expec[f(\rho_{t})]{(f(\rho_1)-f(\rho_0))\wedge C_9}\indb{f\left(\frac{r^{(\epsilon)}}{2}\right)\leq f(d_{t})\leq K}\nonumber\\
=:&\indb{f\left(\frac{r^{(\epsilon)}}{2}\right)\leq f(d_{t})\leq K}c_{7}>0.\nonumber
\end{align}
The last inequality follows from the continuity and positivity ($g(r)>0$ for $r\geq 0$) of the mapping  $r\mapsto \expec[f(\rho_0)=r]{(f(\rho_1)-f(\rho_0))\wedge C_9}$.  \eqref{eq:driftschranke} is now an easy consequence of the following proposition.
\begin{proposition}\label{be:driftschranke2}\mbox{}\\
There is $K\in\mathbb{N}$ such that
$\ce{\left((f(d_{t+1})-f(d_{t}))\wedge C_9\right)\indb{ f(d_{t})> K}}{\mathcal{F}_{t}}\geq\frac{c_6}{2}\indb{ f(d_{t})> K}$.
\end{proposition}
Proof of Proposition~\ref{be:driftschranke2}: Consider for $F\in\mathcal{F}_{t}$ and $K\in\mathbb{N}$ that
\begin{align}\label{eq:grosswachs}
&\expec{\left((f(d_{t+1})-f(d_{t}))\wedge C_9\right)\indb{ f(d_{t})> K}\ind{F}}\nonumber\\
=&\expec{\left((f(d_{t+1})-f(d_{t}))\wedge C_9\right)\ind{\left\{f(d_{t})> K\right\}\cap F\cap\left\{f(d_{t+1})\geq 0 \right\}}}\nonumber\\
&+\expec{\left((f(d_{t+1})-f(d_{t}))\wedge C_9\right)\ind{\left\{f(d_{t})> K\right\}\cap F\cap\left\{f(d_{t+1})<0 \right\}}}\nonumber\\
\geq&\expec{\left((c_{4}(d_{t+1}-d_{t}))\wedge C_9\right)\ind{\left\{f(d_{t})> K\right\}\cap F\cap\left\{f(d_{t+1})\geq 0 \right\}}}\nonumber\\
&+\expec{\left(f\left(\rho_{t+1}\right)-f\left(\rho_{t}\right)\right)\indb{f\left(\rho_{t+1}\right)-f\left(\rho_{t}\right)<-K}\ind{\left\{f\left(d_{t}\right)> K\right\}\cap F\cap\left\{f\left(d_{t+1}\right)<0\right\}}}\nonumber\\
\geq&\expec{\left((c_{4}(d_{t+1}-d_{t}))\wedge C_9\right)\ind{\left\{f(d_{t})> K\right\}\cap F}}\nonumber\\
&+\expec{\expec[f\left(\rho_{t}\right)]{\left(f\left(\rho_{1}\right)-f\left(\rho_{0}\right)\right)\indb{f\left(\rho_{1}\right)-f\left(\rho_{0}\right)<-K}}\ind{\left\{f\left(d_{t}\right)> K\right\}\cap F}}\nonumber\\
\geq&\expec{ \ce{\left(c_4(d_{t+1}-d_{t})\right)\wedge C_9}{\mathcal{F}_{t}} \ind{\left\{f(d_{t})> K\right\}\cap F}}\nonumber\\
&-\expec{\ind{\{f(d_{t})>K\}\cap F}\sum_{n=K}^{\infty}n\prob[f\left(\rho_{t}\right)] {f\left(\rho_{1}\right)-f\left(\rho_{0}\right)<1-n}}\\
=:&\expec{ \ce{\left(c_4(d_{t+1}-d_{t})\right)\wedge C_9}{\mathcal{F}_{t}} \ind{\left\{f(d_{t})> K\right\}\cap F}}-I\nonumber\\
\geq& c_{6}\prob{f(d_{t})>K;F;f(d_{t+1})\geq 0 }-I\nonumber.
\end{align}
For the estimation of $I$ (defined in the above computation) the next lemma ist useful.
\begin{lemma}\label{le:schrumpftail}
For $x,y\in\mathbb{R}^{2}$ and $f\left(\rho^{xy}_{t}\right)\in\mathbb{R}$ we have that\\
\mbox{}\hspace{2cm}$\prob[f(\rho^{xy}_{t})]{f(\rho^{xy}_{0})-f(\rho^{xy}_{1})>n}\leq\frac{\left\|\tilde{g}\right\|_{\infty}}{n\sqrt{2\pi}}\exp\left\{-\frac{n^{2}}{2\left\|\tilde{g}\right\|_{\infty}^{2}}\right\}$.
\end{lemma}
Proof of Lemma~\ref{le:schrumpftail}:
Due $df(\rho^{xy}_{s})=\tilde{g}(\rho^{xy}_{s})dW_{s}+g(\rho^{xy}_{s})ds$ with $\tilde{g}\leq\bar{M}$ and $g\geq0 $ some standard results (~\cite[Proposition 5.2.18]{ks}, \cite[Theorem 4.6]{ks} and \cite[Chapter V, Theorem 1.7]{ry} e.g.) imply
\begin{align}  
\prob[f(\rho^{xy}_{i-1})]{f(\rho^{xy}_{0})-f(\rho^{xy}_{1})>n}\leq&\prob[f(\rho^{xy}_{i-1})]{\int_{0}^{1}\tilde{g}(\rho^{xy}_{s})dW_{s}<-n}\nonumber\\
\leq\prob{W_{\int_{0}^{1}\left\|\tilde{g}\right\|_{\infty}^{2}ds}<-n}\leq\prob{\left\|\tilde{g}\right\|_{\infty} W_{1}<-n}\leq&\frac{\left\|\tilde{g}\right\|_{\infty}}{n\sqrt{2\pi}}\exp\left\{-\frac{n^{2}}{2\left\|\tilde{g}\right\|_{\infty}^{2}}\right\}.\nonumber\end{align} completing the proof of Lemma~\ref{le:schrumpftail}. \hfill$\Box$\\
The convergence 
$I\leq\expec{\ind{\{f(d_{t})>K\}\cap F}\sum_{n=K}^{\infty}\frac{n\left\|\tilde{g}\right\|_{\infty}}{(n-1)\sqrt{2\pi}} \exp\left\{-\frac{(n-1)^{2}}{2\left\|\tilde{g}\right\|_{\infty}^{2}}\right\}} 
\stackrel{K\rightarrow\infty}{\rightarrow}0$ 
implies for sufficiently large $K$ (uniformly in $F$) that\\ $I\leq\frac{c_{4}c_{6}}{2}\prob{\left\{f(d_{t})>K\right\}\cap F\cap\left\{f(d_{t+1})\geq 0 \right\}}$, because we also have that\\ $\cp{f(d_{t+1})\geq0}{f(d_{t})>K}\rightarrow1$ for $K\rightarrow\infty$
which together with (\ref{eq:grosswachs}) completes the proof of Proposition~\ref{be:driftschranke2}.\hfill$\Box$\\
The proof of \eqref{eq:driftschranke} is now straightforward. Choose $K>1$ for which Proposition~\ref{be:driftschranke2} holds,  $c_{7}=c_{7}(K)$ according to \eqref{eq:c7def} and consider
\begin{align}
&\ce{(f(d_{t+1})-f(d_{t}))\wedge C_9}{\mathcal{F}_{t}}\nonumber\\
=&\ce{\left((f(d_{t+1})-f(d_{t}))\wedge C_9\right)\indb{f(d_{t})<f\left(\frac{r^{(\epsilon)}}{2}\right)}}{\mathcal{F\rule{0mm}{4mm}}_{t}}\nonumber\\
&+\ce{ \left((f(d_{t+1})-f(d_{t}))\wedge C_9\right) \indb{ f\left(\frac{r^{(\epsilon)}}{2}\right) \leq f(d_{t})\leq K} }{\mathcal{F\rule{0mm}{4mm}}_{t}}\nonumber\\
&+\ce{\left((f(d_{t+1})-f(d_{t}))\wedge C_9\right)\indb{f(d_{t})> K}}{\mathcal{F}_{t}}\nonumber\\
\geq&\left[f\left(r^{(\epsilon)}\right)-f\left(\frac{r^{(\epsilon)}}{2}\right)\right]\wedge\frac{1}{8}(\beta_{N}(d-1)-\beta_{L})\wedge C_9\wedge c_{7}\wedge\frac{c_6}{2}=:C_8>0 \nonumber.
\end{align}
The proof of Lemma~\ref{le:exf} is complete. \hfill$\Box$\\
The estimate (on average) is about to be transformed into one of the probability of the event that our original set is not large after a long time.
\subsubsection{Pathwise Growth Of $f(d_{t})$}\label{suse:pww}
We have $\ce{f(d_{t+1})\wedge (f(d_{t})+C_9)}{\mathcal{F}_{t}}-f(d_{t})\geq C_8>0$. So we can verify the assumptions of Lemma~\ref{le:mmb} for $\xi_{i}:=f(d_{i-1})-[f(d_{i})\wedge (f(d_{i-1})+C_9)]+C_8$. we only have to prove $\mathbb{E}[\xi_{i}^{m}]\leq K_{m}$ for certain real $K_{m}$.
\begin{align}
\expec{|\xi_{i}|^{m}}\leq&2^{m}C_8^m+ 2^mC_9^m\prob{d_{i}>d_{i-1}}+2^{m}\expec{\left( f(d_{i-1})-f(d_{i})\right)^{m} \indb{d_{i-1}>d_{i}}}\nonumber\\
\leq&2^{m}C_8^m+ 2^mC_9^m+2^m\expec{\left( f(d_{i-1})-f(d_{i})\right)^{m} \indb{d_{i-1}>d_{i}}}\nonumber\\
:=&2^m(C_8^m+C_9^m)+2^m\mathit{I}\label{eq:ximoment1}.
\end{align}
To estimate $I$ we choose $x$ and $y$ in $\gamma$ such that $\left\|x_{i-1}-y_{i-1}\right\|=d_{i-1}$. For shrinking of $d_{t}$ implies decreasing of  $||x_{t}-y_{t}||$ we can further conclude
\begin{align}  
I\leq&\expec{\left( f(\rho^{xy}_{i-1})-f(\rho^{xy}_{i})\right)^{m} \indb{\rho^{xy}_{i-1}>\rho^{xy}_{i}}}\nonumber\\
=&\expec{\ce{ \left( f(\rho^{xy}_{i-1})-f(\rho^{xy}_{i})\right)^{m} \indb{f\left(\rho^{xy}_{i-1}\right)>f\left(\rho^{xy}_{i}\right)}}{\rule{0mm}{4mm}\mathcal{F}_{i-1}} }\nonumber\\
\leq&1+\sum_{n=1}^{\infty}(n+1)^{m}\sup_{f(r)\in\mathbb{R}^{+}}\prob[f(r)]{f(\rho^{xy}_{0})-f(\rho^{xy}_{1})>n}\label{eq:ximoment2}.
\end{align}  
Combining (\ref{eq:ximoment1}), (\ref{eq:ximoment2}) and Lemma~\ref{le:schrumpftail} yields
\[
\expec{|\xi_{i}|^{m}}\leq2^{m}(C^m_8+C_9^m+1)+2^{m}\sum_{n=1}^{\infty}\frac{\left\|\tilde{g}\right\|_{\infty}(n+1)^{m}}{n\sqrt{2\pi}}\exp\left\{-\frac{n^{2}}{2\left\|\tilde{g}\right\|_{\infty}^{2}}\right\}=:K_{m}.
\]
Concluding with Lemma~\ref{le:mmb} we have for $m\in\mathbb{N}$ the existence of $\kappa^{(1)}_{m}\in\mathbb{R}$, such that for $n\geq\frac{2}{C_8}\left|f(d_{0})\right|$ the following holds.
\begin{align}  
\prob{d_{n}<1}=&\prob{f(d_{n})<0}=\prob{\sum_{i=0}^{n-1}\left(f(d_{i})-f(d_{i+1})+C_8\rule{0pt}{12pt}\right)> f(d_{0})+C_8  n}\nonumber\\
\leq&\prob{\sum_{i=0}^{n-1}\left(f(d_{i})-[f(d_{i+1})\wedge (f(d_{i})+K)]+C_8 \rule{0pt}{12pt} \right)\geq f(d_{0})+C_8 n}\nonumber\\
=&\prob{\sum_{i=1}^{n}\xi_{i}\geq f(d_{0})+C_8 n}\leq\prob{\sum_{i=1}^{n}\xi_{i}\geq \frac{C_8 n}{2}}\leq\kappa_{m}^{(1)} n^{-m}.
\end{align}
Increasing $\kappa_{m}^{(1)}$ ensures that this holds for all $n$ (this is to be assumed). \\ 
\textbf{Remark:} The assumption of largeness of $\gamma$ makes this correction uniform in $\gamma$.\\
So we can estimate the probability of
$F_{n}:=\{\exists i\in[\left\lfloor \sqrt{n}\right\rfloor,\infty]\cap\mathbb{N}:d_{i}<1\}$ via
\begin{equation}\label{eqn:kappaquer2komma5}
\prob{F_{n}}\leq\sum_{i=\left\lfloor \sqrt{n}\right\rfloor}^{\infty}\kappa^{(1)}_{2+2m}i^{-2-2m}
\leq\left(\kappa^{(1)}_{2+2m}\sum_{i=1}^{\infty}i^{-2}\right)n^{-m}=:\kappa^{(6)}n^{-m}.
\end{equation}
A simple Borel-Cantelli-argument shows that the flow cannot contract a non-rivial set to a point i.e. $d_t$ a.s. does not converge to zero as $t\rightarrow\infty$.
\subsubsection{Getting Estimates On The Tails Of $\tau^{R}(\gamma,P)$}
First let $r_t:=\dist{(\gamma_t,P)}$ and observe for $n\in\mathbb{N}$
\be  \prob{\tau^{R}(\gamma,P)>n}\leq\prob{F_{n}}+\mathbb{P}[F_{n}^{C},\bigcap_{i=\left\lfloor \sqrt{n} \right\rfloor}^{n}\left\{r_{i}>R\right\}]=:I+\mathit{II}.\label{eqn:kappaquer1}\ee 
$I$ is aready treated, so only $\mathit{II}$ is left. For arbitrary $\delta>0$ and $n\geq4\vee4(r_{0}-R)\delta^{-1}$ we can estimate
{\fontsize{10pt}{2ex}\[
\hspace*{-2cm}\mathit{II}\leq\prob{ \bigcap_{i=\left\lfloor n \right\rfloor}^{n}\left\{r_{i}>R\right\},\bigcap_{i=\left\lfloor\sqrt{n}\right\rfloor}^{\infty}\left\{d_{i}\geq1\right\} } \hspace{5.17cm}\]\[
=\prob{(r_{\left\lfloor \sqrt{n}\right\rfloor}-r_{0})+\sum_{i=\left\lfloor \sqrt{n}+1\right\rfloor}^{n}\left(r_{i}-r_{i-1}\right)>R-r_{0},\bigcap_{i=\left\lfloor n \right\rfloor}^{n-1}\left\{r_{i}>R\right\},\bigcap_{i=\left\lfloor\sqrt{n}\right\rfloor}^{\infty}\left\{d_{i}\geq1\right\}} \hspace{0.99cm}\]\[
\leq\prob{\eta^{(n)}\geq\frac{\delta n}{8}}\hspace{11.2cm}\]\[
+\prob{\sum_{i=\left\lfloor \sqrt{n}+1\right\rfloor}^{n}\left[(r_{i}-r_{i-1})\indb{d_{i-1}\geq1,r_{i-1}>R}-\delta\ind{\left\{d_{i-1}<1\right\}\cup\left\{r_{i-1}\leq R\right\}}+\delta\right]\geq\frac{\delta}{4}(n-\left\lfloor\sqrt{n}\right\rfloor)}\]\begin{equation}
=:\mathit{III}+\mathit{IV}. \hspace{11.4cm}\label{eqn:kappaquer2}
\end{equation}}
Therein  $\eta^{(n)}:=r_{\left\lfloor \sqrt{n}\right\rfloor}-r_{0}$ is used.
The term $\mathit{III}$ can be estimated by the growth of a Brownian motion. Choose $z\in\gamma$ with $\left\|z-P\right\|=r_{0}$. Then we have
\be 
\mathit{III}\leq\prob{\left\| z_{\left\lfloor \sqrt{n}\right\rfloor}-z_{0} \right\|\geq\frac{\delta n}{8}}\leq\kappa^{(3)}_{m}n^{-m}\label{eqn:kappaquer3}
\ee 
for suitable $\kappa^{(3)}_{m}\in \mathbb{R}$.
The estimation of $\mathit{IV}$ applies Lemma~\ref{le:mmb} again. For $\delta>0$ and $n\geq4\vee4(r_{0}-R)\delta^{-1}$ observe
\be \label{eq:badsum}
\mathit{IV}\leq\prob{\sum_{i=\left\lfloor \sqrt{n}+1\right\rfloor}^{n}\xi_{i}^{(C_{10},\delta)} \geq(n-\left\lfloor \sqrt{n}\right\rfloor)\frac{\delta}{4}}.
\ee
Therein for $C_{10}>0$ and $i\in\mathbb{N}$ set\\
$\xi_{i}^{(C_{10},\delta)}:=\left[\left(r_{i}\vee (r_{i-1}-C_{10})-r_{i-1}\right)\indb{d_{i-1}\geq1,r_{i-1}>R}-\delta\ind{\left\{d_{i-1}<1\right\}\cup\left\{r_{i-1}\leq R\right\}}+\delta\right]$
The sequel aims at showing that $(\xi_{i}^{(C_{10},\delta)}:i\in\mathbb{N})$ for suitable $C_{10}$ and $\delta$ satisfies the assumptions of Lemma~\ref{le:mmb}. Afterwards this lemma and a treatment of the fact that there are some terms $\xi_{i}^{(C_{10},\delta)}$ missing in the last sum which makes Lemma~\ref{le:mmb} not directly suitable for (\ref{eq:badsum}) will complete the proof.
Therefore we have to show: $\expec{\left|\xi_{i}^{(C_{10},\delta)}\right|^{m}}\leq K_{m}<\infty$ for any $m$ and uniformly in $i$. 
\be  
\expec{\left|\xi_{i}^{(C_{10},\delta)}\right|^{m}}\leq2^m\delta^m+2^mC_{10}^m+2^m\expec{\left(r_{i}-r_{i-1}\right)^{m}\indb{d_{i-1}\geq1,r_{i}>r_{i-1}}}\nonumber.
\ee  
For $\gamma$ cannot get away from $P$ without having its nearest (w.r.t. $P$) point doing so we can proceed for the estimation of the above as follows: Let $z\in\gamma$ such that $\left\|z_{i-1}-P\right\|=r_{i-1}$ and consider
\be  
\expec{\left(r_{i}-r_{i-1}\right)^{m}\indb{d_{i-1}\geq1,r_{i}>r_{i-1}}}\leq\expec{|z_i-z_{i-1}|^m} =\expec{\left|\mathcal{N}(0,1)^{\otimes d}\right|^{m}}.
\ee  
Therefore we can choose
$K_{m}=K_{m}(C_{10},\delta):=2^{m}\delta^{m}+2C_{10}^{m}+2^{m}\expec{\left|\mathcal{N}(0,1)^{\otimes d}\right|^{m}}$ and it only remains to show that there are $C_{10}>0$ and $\delta>0$ such that $\ce{\xi_{i}^{(C_{10},\delta)}}{\xi_{i-1}^{(C_{10},\delta)}\ldots\xi_{1}^{(C_{10},\delta)}}$ is negative. ($\expec{\left|\mathcal{N}(0,1)^{\otimes d}\right|^{m}}$ here simply denotes the $m$th moment of the $d$-dimensional normal.)
Therefore it is sufficient to show $\ce{\xi_{i}^{(C_{10},\delta)}}{\rule[2mm]{0mm}{2mm}\mathcal{F}_{i-1}}\leq0$ for suitable $C_{10}$ and $\delta$. On $\left\{d_{i-1}<1\right\}$ and on $\left\{r_{i-1}\leq R\right\}$ this is evident. On  $\left\{d_{i-1}\geq1, r_{i-1}>R\right\}$ we use Lemma~\ref{le:lokdrift}. Because of 2. there is $0.5\geq\rho>0$ with $G''(1,\rho)=:2\delta>0$. 1. yields the existence of $C_{10}>0$ and $\hat{r}>0$ such that we have for $r>\hat{r}$: $G'(C_{10},1,r)<\delta$. Now choose $x,y\in\gamma$ with $\left\|x_{i-1}-P\right\|=r_{i-1}$ and $\left\|y_{i-1}-x_{i-1}\right\|=\rho$. With 3. we conclude ($\tau\equiv i-1$ and $z\equiv P$) that for
$l_{1}:=\left\|x_{i-1}-P\right\|\wedge\left\|y_{i-1}-P\right\|$ and $l_{2}:=\left\|x_{i}-P\right\|\wedge\left\|y_{i}-P\right\|$ we have
\begin{align}
&\ce{\left(r_i\vee(r_{i-1}-C_{10})-r_{i-1}\right)\indb{r_{i-1}>R,d_{i-1}\geq1}}{\mathcal{F}_{i-1}}\nonumber\\
\leq&\ce{\left(l_2\vee(l_1-C_{10})-l_1\right)}{\mathcal{F}_{i-1}}\indb{r_{i-1}>R,d_{i-1}\geq1}\nonumber\\
\leq& \left(G'\left(C_{10},1,r_{i-1}\right)-G''(1,\rho)\right)\indb{r_{i-1}>R,d_{i-1}\geq1}
\leq-\delta\indb{r_{i-1}>R,d_{i-1}\geq1},\nonumber
\end{align} 
provided we choose $R:=\hat{r}$ (which we do). So we can apply Lemma~\ref{le:mmb} to $(\xi_{i}^{(\delta,C_{10})}:i\in\mathbb{N})$ for these $C_{10}$ and $\delta$. We will abbreviate $\xi_{i}^{(C_{10},\delta)}$ as $\xi_{i}$. Fix $C_{10}$ and $\delta$ satisfiing the assumptions of Lemma~\ref{le:mmb} and conclude for $n\geq4\vee4(r_{0}-R)\delta^{-1}\vee\left(16 C_{10}+1\right)^{2}\delta^{-2}$:
For $m\in\mathbb{N}$ there is $\kappa^{(4)}_{m}\in\mathbb{R}$ such that 
\begin{align}
\mathit{IV}\leq&\prob{\sum_{i=\left\lfloor \sqrt{n}+1\right\rfloor}^{n}\xi_{i} \geq(n-\left\lfloor\sqrt{n}\right\rfloor)\frac{\delta}{4},\sum_{i=1}^{n}\xi_{i} \geq\frac{\delta n}{16}}\nonumber\\
&+\prob{\sum_{i=\left\lfloor \sqrt{n}+1\right\rfloor}^{n}\xi_{i} \geq(n-\left\lfloor\sqrt{n}\right\rfloor)\frac{\delta}{4},\sum_{i=1}^{\left\lfloor \sqrt{n}\right\rfloor}\xi_{i} \leq-\frac{\delta n}{16}}\nonumber\\
\leq&\prob{\sum_{i=1}^{n}\xi_{i} \geq\frac{\delta n}{16}}+\prob{\sum_{i=1}^{\left\lfloor \sqrt{n}\right\rfloor}\xi_{i} \leq-\frac{\delta n}{16}}\leq\kappa^{(4)}_{m}n^{-m}.\label{eqn:kappaquer4}
\end{align}
Observe that due to $n>\left(16C_{10}\right)^{2}\delta^{-2}\Rightarrow \sum_{i=1}^{\left\lfloor \sqrt{n}\right\rfloor}\frac{\xi_{i}}{\left\lfloor\sqrt{n}\right\rfloor}\geq -C_{10}>-\frac{\delta \sqrt{n}}{16} $ the last term vanishes. Combining the equations (\ref{eqn:kappaquer1}), (\ref{eqn:kappaquer2komma5}), (\ref{eqn:kappaquer2}),  (\ref{eqn:kappaquer3}) and (\ref{eqn:kappaquer4}) yields for $n\geq4\vee4(r_{0}-R)\delta^{-1}\vee\left(16C_{10}+1\right)^{2}\delta^{-2}$:
\be  
\prob{\tau^{R}(\gamma,P)>n}\leq\kappa^{(6)}_{m}n^{-m}+\kappa^{(3)}_{m}n^{-m}+\kappa^{(4)}_{m}n^{-m}=:\kappa_m^{(5)}n^{-m},
\ee  
which proves that for $m\in\mathbb{N}$ the choice\vspace{-1mm}
\[\kappa^{(2)}_{m}:=\left[(\kappa^{(5)}_{m}\vee1)\sup_{r>1}\left(\frac{r}{\left\lfloor r\right\rfloor}\right)^{m}\right]\left[ 4\vee\left(\frac{16C_{10}+1}{\delta}\right)^{2}\vee\frac{4}{\delta}\right]<\infty \]\vspace{-1mm} is appropriate completing the proof of Theorem~\ref{th:cor6}.\hfill$\Box$\newline
\subsection{Linear Expansion And Stable Norm}
The next two subsections follow closely the line of thought of~\cite{dkk} although we cannot use their results directly.
\subsubsection{Implications Of Theorem~\ref{th:cor6}}
For collecting the following corollaries of Theorem~\ref{th:cor6} we let
$$  
\mathcal{W}_{t}(\gamma):=\bigcup_{0\leq s\leq t}\gamma_{s}\textnormal{ and }
\mathcal{W}_{t}^{R}(\gamma):=\left\{ x\in\mathbb{R}^d: \dist\left(x,\mathcal{W}_{t}(\gamma)\right)\leq R\right\}. 
$$
\begin{corollary}\label{ko:wrt}
There are positive constants $C_{11}$ and $R$, such that $\mathbb{P}\textnormal{-}a.s.$ we have for large $t$ (i.e. for all $t$ that are bigger than some a.s. finite random variable) that
\begin{equation}
K_{C_{11}t}(0)\subset\mathcal{W}_{t}^{R}(\gamma).
\end{equation}
$K_{r}(x)$ denotes the closed $r$-Ball centered at $x$ as before. 
\end{corollary}
Proof: Cover $K_{C_{11}t}(0)$ with balls of radius $R$. Due to Theorem~\ref{th:cor6} the probability, that a fixed one of these balls has not been hit by $\gamma$ up to time $t$, decays faster than any power of $t$, if we choose $C_{11}$ small enough and $R$ large enough. For the number of balls needed to cover $K_{C_{11}t}(0)$ only grows like $t^{d}$ the probability that any of these balls has not been hit up to time $t$ decays faster than any power of $t$ provided $R$ is sufficiently large and $C_{11}$ sufficiently small. So Corollary~\ref{ko:wrt} follows from the first Borel-Cantelli-lemma.\hfill$\Box$\\
For the sequel fix $R>0$ large enough for Theorem~\ref{th:cor6} and Corollary~\ref{ko:wrt} to hold with this $R$. Assuming that $\gamma$ is large makes all the estimates of Theorem~\ref{th:cor6} uniform in $\gamma\in\mathcal{C}_{R}$ with $\mathcal{C}_{R}:=\left\{\gamma: \diam(\gamma)\geq1 ,\gamma\subset K_{2R}(0)\right\}$. (w.l.o.g. we assume $R>1$ ). The following is immediate from Theorem~\ref{th:cor6}.
\begin{corollary}\label{ko:ggi}
The family of random variables $\left(\left(\frac{\tau^{R}(\gamma, tv)}{t}\right)^{k}\right)_{t\geq1,\left\|v\right\|=1,\gamma\in\mathcal{C}_{R}}$ is  uniformly integrable for any $k\in\mathbb{N}$. 
\end{corollary}
\subsubsection{The Stable Norm}
\begin{figure}\centering \resizebox{6cm}{!}{\includegraphics{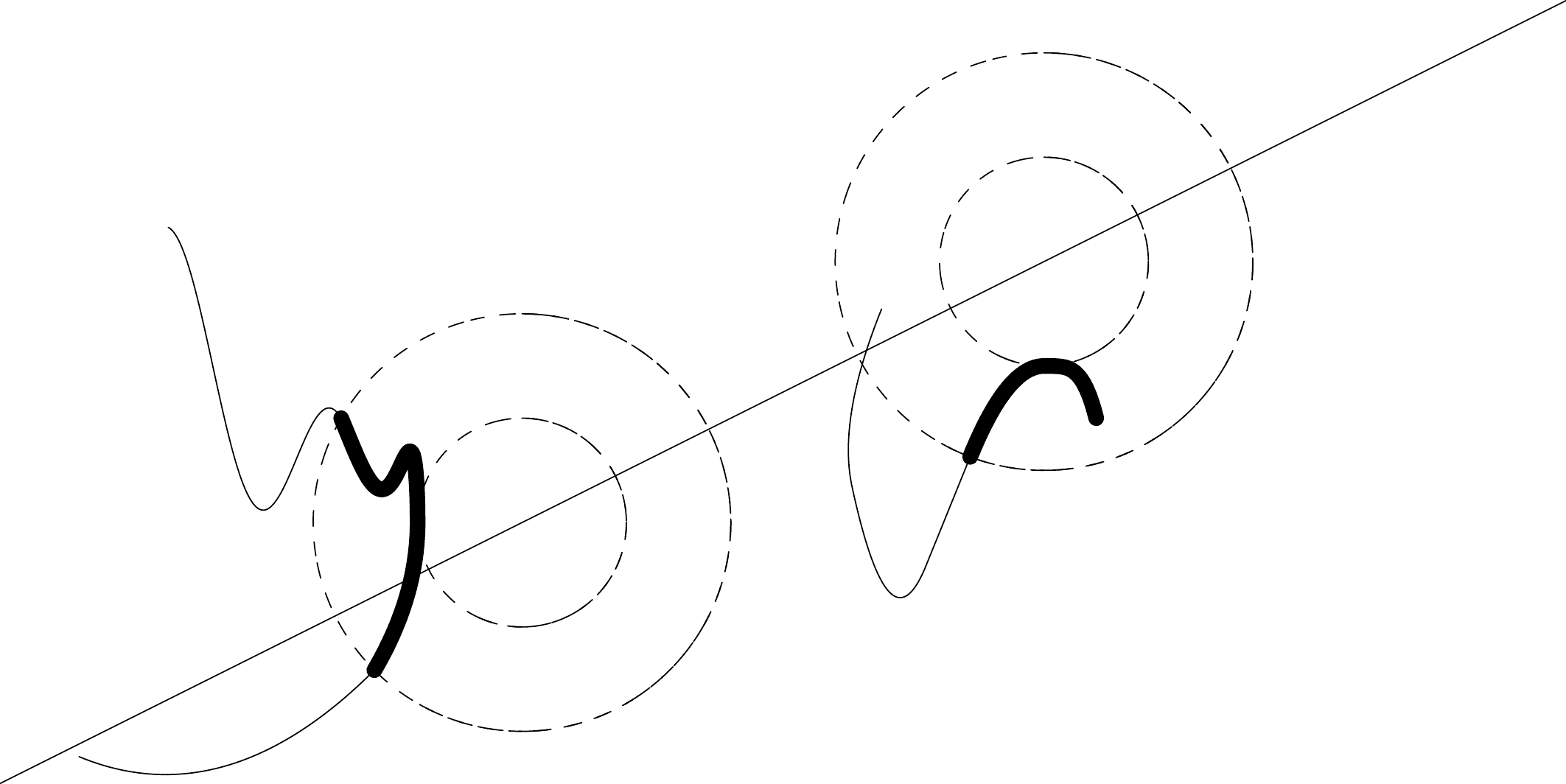}}\caption{line: direction of $v$, fat: $\gamma^{(i)}$, regular: rest of $\Phi_{\tau_{i}^{R}}\gamma$\label{abb:uni}}\end{figure}

Set $\left|v\right|^{R}:=\sup_{\gamma\in\mathcal{C}_{R}}\expec{\tau^{R}(\gamma,v)}$, which due to the isotropic properties of the flow does not depend on the direction of $v$. We obviously have
\begin{equation}\label{eq:suba1}
\expec{\tau^{2R}\left(\gamma,(t_{1}+t_{2})v\right)}\leq \expec{\tau^{R}\left(\gamma,t_{1}v\right)}+\sup_{\check{\gamma}\in\mathcal{C}_{R}}\expec{\tau^{R}\left(\check{\gamma},t_{2}v\right)}.
\end{equation} With Theorem~\ref{th:cor6} we get in addition
\begin{equation}\label{eq:suba2}
\expec{\tau^{R}\left(\gamma,(t_{1}+t_{2})v\right)}\leq\expec{\tau^{2R}\left(\gamma,(t_{1}+t_{2})v\right)}+C_{12}
\end{equation}
for some constant $C_{12}>0$. Combining (\ref{eq:suba1}) and (\ref{eq:suba2}) yields the subadditivity of $t\mapsto \left|tv\right|^{R}+C_{12}$. Using Feketes lemma we conclude that
$\left\|v\right\|^{R}:=\lim_{t\rightarrow\infty}(\left|tv\right|^{R}+C_{12})t^{-1}=\lim_{t\rightarrow\infty}\left|tv\right|^{R}t^{-1}$ is well-defined i.e. the limit exists and equals $\inf_{t\geq0} (\left|tv\right|^{R}+C_{12})t^{-1}$.
Since $|v|^{R}$ only depends on $\left\|v\right\|$ and since it is increasing with respect to this argument we get (again from the isotropy of the flow) that
$\left\|s v_{1}+(1-s)v_{2}\right\|^{R}\leq s\left\|v_{1}\right\|^{R}+(1-s)\left\|v_{2}\right\|^{R}$. Set $\mathcal{B}:=\{v\in\mathbb{R}^{d}:\left\|v\right\|^{R}\leq1\}$ and observe that $\mathcal{B}$ is a compact convex set (see Lemma~\ref{le:subgauss}). Corollary~\ref{ko:wrt} shows $\left\|v\right\|^{R}\neq0$ provided $v\neq0$. Of course the isotropic properties of the flow imply that $\mathcal{B}$ is a ball centered at the origin. We will show later that its radius does not depend on $R$.
First we can prove the following lemma.
\begin{lemma}\label{le:usvor}
For any $\gamma\in\mathcal{C}_{R}$ and $\epsilon>0$ there is $\mathbb{P}$-a.s. $T(\gamma,\epsilon)>0$, such that for $t>T(\gamma,\epsilon)$ we have:
$(1-\epsilon)t\mathcal{B}\subset\mathcal{W}^{R}_{t}(\gamma)$.
\end{lemma}
Proof: We need to show that for $v$ with $\left\|v\right\|^{R}\leq1$ and $m\in\mathbb{N}$ there is $\kappa_m^{(7)}=\kappa_m^{(7)}(\epsilon)>0$, such that  \begin{equation}\label{eq:hittail}\prob{\tau^{R}(\gamma,tv)\geq(1+\epsilon)t}\leq \kappa_m^{(7)}t^{-m}\end{equation} holds uniformly in $\gamma\in\mathcal{C}_{R}$ and $\left\|v\right\|^{R}\leq1$.  All the estimates we made so far are uniform in $\left\|v\right\|^{R}=1$, because they do not depend on the direction of $v$.
By definition of $\left\|.\right\|^{R}$ there is $\tilde{t}>0$ with $\expec{\tau^{R}(\gamma,tv)}\leq(1+\frac{\epsilon}{2})t$ for any $t\geq\tilde{t}$ and $\gamma\in\mathcal{C}_{R}$. Define the stopping time $\tau_{1}^{R}$ via \[\tau^{R}_{1}:=\inf\left\{t>0,\gamma_{t}\cap K_{R}(\tilde{t}v)\neq\emptyset,\diam(\gamma_{t})\geq1\right\}.\] Denote by $\gamma^{(1)}$ a large connected subset of  $\gamma_{\tau_{1}^{R}}$ which is contained in $K_{2R}(\tilde{t}v)$ and which has non-empty intersection with $K_{R}(\tilde{t}v)$. We can choose it to be $\mathcal{F}_{\tau^{R}_{1}}$-mesureable which we do. Now define an increasing sequence of stopping times $\left(\tau^{(i)}:i\in\mathbb{N}\right)$ recursively via
{\fontsize{10pt}{2ex} 
\begin{align}  \label{eq:tautimes} \tau_{i}^{R}:=&\inf\left\{t>\tau_{i-1}^{R},\diam\left(\Phi_{\tau_{i-1}^{R},t}\left(\gamma^{(i-1)}\right)\right) \geq1, \Phi_{\tau_{i-1}^{R},t}\left(\gamma^{(i-1)}\right)\cap K_{R}(i\tilde{t}v)\neq\emptyset\right\}, \nonumber\\
\gamma^{(i)}:=&\Phi_{\tau_{i-1}^{R},\tau_{i}^{R}}\left(\gamma^{(i-1)}\right)\cap K_{2R}(\tilde{t}iv).\nonumber
\end{align}  }
\hspace{-1.5mm}(If necessary we choose a subset of $\gamma^{(i)}$ as $\gamma^{(i)}$ to ensure that it is connected.) We have (putting $\tau_{0}^{R}\equiv0$) that $\tau^{R}(\gamma,n\tilde{t}v)\leq\sum_{j=1}^{n}(\tau_{j}^{R}-\tau_{j-1}^{R})$ (see Fig.~\ref{abb:uni}).
Due to the strong Markov-property, the isotropy of $\Phi$ and the definition of $\tilde{t}$ yield $\ce{\tau_{j}^{R}-\tau_{j-1}^{R}}{\mathcal{F}_{\tau_{j-1}^{R}}}\leq\left(1+\frac{\epsilon}{2}\right)\tilde{t}$. Due to Theorem~\ref{th:cor6} we can define $\xi_{j}:=\tau_{j}^{R}-\tau_{j-1}^{R}-\left(1+\frac{\epsilon}{2}\right)\tilde{t}$ and obtain that the sequence $\left(\xi_{i}:i\in\mathbb{N}\right)$ satisfies the assumptions of Lemma~\ref{le:mmb}. So we conclude
\begin{align}
\prob{\tau^{R}(\gamma,n\tilde{t}v)\geq(1+\epsilon)n\tilde{t}}\leq&\prob{\sum_{j=1}^{n}\tau_{j}^{R}-\tau_{j-1}^{R}\geq(1+\epsilon)n\tilde{t}}\nonumber\\
=\prob{\sum_{j=1}^{n}\xi_{j}\geq\frac{\epsilon}{2}n\tilde{t}}
\leq\kappa_m^{(1)}n^{-m}=&\kappa_m^{(1)}\tilde{t}^{m}(n\tilde{t})^{-m}
=:\kappa_m^{(7)}\left(n\tilde{t}\right)^{-m}
\end{align} 
which implies that a.s. for any $\epsilon>0$ the inclusion $(1-\epsilon)n\tilde{t}\mathcal{B}\subset\mathcal{W}^{R}_{n\tilde{t}}(\gamma)$ fails to hold only a finite number of times. Considering that the definition
$t^{\downarrow}:=\left\lfloor\frac{t}{\tilde{t}}\right\rfloor\tilde{t}$ implies $
\lim_{t\rightarrow\infty}\frac{t^{\downarrow}}{t}=1$ and for all $\epsilon>0$ that there is $\check{t}>0$ such that for $t\geq\check{t}$ we have $t^{\downarrow}\geq\frac{1-\epsilon}{1-\frac{1}{2}\epsilon}t$
 finally proofs for $t\geq\check{t}\vee \max\left\{n\in\mathbb{N}:(1-\frac{\epsilon}{2})n\tilde{t}\mathcal{B} \nsubseteq \mathcal{W}^{R}_{n\tilde{t}}(\gamma)\right\}\tilde{t}$ that $(1-\epsilon)t\mathcal{B}\subset\left(1-\frac{\epsilon}{2}\right)t^{\downarrow}\mathcal{B}\subset\mathcal{W}^R_{t^{\downarrow}}\subset\mathcal{W}^R_t\textnormal{ a.s.}$ and hence completes the proof of Lemma~\ref{le:usvor}.\hfill$\Box$\\
\subsection{Sweeping Lemma And Lower Bound - The $2$-Dimensional Case}
In this subsection assume $d=2$. We will also assume that $\gamma$ is a curve (which we could have assumed before). In this case we have
\begin{theorem}\label{sa:us}
For any $\gamma\in\mathcal{C}_{R}$ and $\epsilon>0$ there is $\mathbb{P}$-a.s. $T(\gamma,\epsilon)>0$, such that for any $t>T(\gamma,\epsilon)$ the following holds
$$(1-\epsilon)t\mathcal{B}\subset\mathcal{W}_{t}(\gamma).$$
This is 1. of Theorem\ref{th:main}.
\end{theorem} 
(Note that we do not distinguish between the $T(\gamma,\epsilon)$ here and the $T(\gamma,\epsilon)$ of Lemma~\ref{le:usvor} because the two times are very close to each other as we will see in the sequel.
The proof of Theorem~\ref{sa:us} depends apart from Lemma~\ref{le:usvor} on the following Sweeping Lemma, which will be proved after Theorem~\ref{sa:us}. 
\begin{lemma}\label{le:sweep}
Let $\gamma$ be a large curve with $\dist(\gamma,P)\leq R$ (for an $R$ as defined before). Define $\tilde{\tau}=\tilde{\tau}^{R}(\gamma,P):=\tilde{\tau}(P):=\inf_{t>0}\{K_{R}(P)\subset\bigcup_{0\leq s\leq t}\gamma_{s}\}$
Then for $m\in\mathbb{N}$ there is $\kappa^{(8)}_m\in\mathbb{R}$ such that 
$\prob{\tilde{\tau}>t}\leq \kappa_m^{(8)}t^{-m}$ holds uniformly in $\gamma$.\end{lemma}
Proof of Theorem~\ref{sa:us}:  
There is  a positive integer $k$, such that for any $n\in\mathbb{N}$ $(1-\epsilon)n\mathcal{B}$ can be covered with $n^{2}k$ balls $\left\{K_{R}\left(P_{i}^{n}\right):i=1,\ldots,n^{2}k\right\}$ of radius $R$. By~\eqref{eq:hittail} the probability, that one of these balls has not been hit by the (at the hitting time large) curve $\gamma$ up to time $(1-0,5\epsilon)n$ decays faster than any power of $n$. Due to Lemma~\ref{le:sweep} $\prob{\tilde{\tau}\left(P_i^n\right)-\tau_{R}\left(\gamma,P_i^n\right)\geq 0,5\epsilon n}$ decays faster than any power of $n$, too. So the probability, that there is one among the balls $K_{R}\left(P_{i}^{n}\right)$ for $i\in\left\{1,\ldots,n^{2}k\right\}$ that is not completely included in $\mathcal{W}_{n}$ at time $n$ decays faster than any power of $t$, which proves Theorem~\ref{sa:us} because we have for large $t$ that 
$(1-2\epsilon)t\mathcal{B}\subset(1-\epsilon)\left\lfloor t\right\rfloor\mathcal{B}\subset\mathcal{W}_{\left\lfloor t\right\rfloor}\subset\mathcal{W}_{t}$.\hfill$\Box$\\
Proof of Lemma~\ref{le:sweep}: The proof consists of six steps. These are carried out similarly to a proof in \cite{dkk}. 
\subsubsection{Localizing Of Lemma~\ref{le:sweep}}
Assume we can prove the following: For any $Q\in K_{R}(P)$ there is an open superset $U_{Q}$ of $Q$, such that for any $\tilde{\tau}_{Q}:=\inf_{t>0}\left\{U_{Q}\subset\cup_{0\leq s\leq t}\gamma_{s}\right\}$ the following holds. For $m\in\mathbb{N}$ there is $\kappa_m^{(9)}\in\mathbb{R}$, such that 
\begin{equation}\label{eq:lokal}\prob{\tilde{\tau}_{Q}>t}\leq \kappa_m^{(9)} t^{-m}\end{equation} holds uniformly for large curves $\gamma$ which have an non-empty intersection with $K_{R}(P)$.
For the covering of $\overline{K_{R}(P)}$ requires only a finite number of the $U_{Q}$ Lemma~\ref{le:sweep} holds because of $\left\{\tilde{\tau}>t\right\}\subset \left\{\tilde{\tau}_{Q}>t \textnormal{ for one of these } Q \right\}$. 
\subsubsection{Definition Of A Small Square}
Set $(q_{1},q_{2}):=Q\in K_{R}(P)$ and consider the following elements of the RKHS $\mathcal{H}$ of $\Phi$:
\bee  
V_{1}^{i}(.)&:=&\int b^{ij}(.-y)d\delta_{Q}\otimes\delta_{1}(y,j)=b^{i,1}(.-Q):i=1,2;\nonumber\\
V_{2}^{i}(.)&:=&\int b^{ij}(.-y)d\delta_{Q}\otimes\delta_{2}(y,j)=b^{i,2}(.-Q):i=1,2.
\eee  
We have
$w_{1}:=V_{1}(Q)=b^{.1}(0)=\left(\begin{array}{c}1\\0\end{array}\right)$  and $w_{2}:=V_{2}(Q)=b^{.2}(0)=\left(\begin{array}{c}0\\1\end{array}\right)$.
Lemma~\ref{le:isoko} implies the Taylor expansions \\
$\left(V_{1}(x-Q),V_{2}(x-Q)\right)=E_{2}+O\left(\left\|x-Q\right\|^{2}\right):(x-Q\rightarrow0)$. So there are $C_{13}>0$ and $\delta>0$ such that we have $\left\|V_{1/2}(x)-w_{1/2}\right\|\leq C_{13} \left\|x-Q\right\|^{2}$ for $\left\|x-Q\right\|<\delta$.
This implies that for $n\in\mathbb{N}$ there is $\epsilon>0$ such that \[U_{Q}^{n,\epsilon}:=\left]q_{1}-n\epsilon,q_{1}+n\epsilon\right[\times\left]q_{2}-n\epsilon,q_{2}+n\epsilon\right[\subset\left\{y\in\mathbb{R}^{2}:\left\|V_{1/2}(y)-w_{1/2}\right\|\leq\epsilon\right\},\]
because for $\epsilon\leq2^{-0,5}n^{-1}\delta\wedge\left(2C_{13}n^{2}\right)^{-1}$ and $x:=(x_{1},x_{2})\in U_{Q}^{n,\epsilon}$ we get 
\[
\left\|V_{1}(x)-w_{1}\right\|\vee\left\|V_{2}(x)-w_{2}\right\|\leq C_{13}\left[(x_{1}-q_{1})^{2}+(x_{2}-q_{2})^{2}\right]\leq 2C_{13}n^{2}\epsilon^{2}\leq\epsilon.
\]  
Note that this still holds, if we decrease $\epsilon$ (for a fixed $n$). Define
\begin{align}  
\tilde{U}_{Q}^{n}:=&\left]q_{1}-\frac{n\epsilon}{2},q_{1}+\frac{n\epsilon}{2}\right[\times\left]q_{2}-\frac{n\epsilon}{2},q_{2}+\frac{n\epsilon}{2}\right[,\nonumber\\
t_{u}^{n}:=&\frac{n\epsilon}{2}\left(\sup_{z\in U_{Q}^{n,\epsilon}}\left(\left\|V_{1}(z)\right\|\vee\left\|V_{2}(z)\right\|\right)\right)^{-1}>0.\nonumber
\end{align}
We may assume $t_{u}^{n}\geq3^{-1}n\epsilon$ as well as $\epsilon\leq 102^{-1}$ (otherwise choose a smaller $\epsilon$). Denote by $\psi^{(i)}_{st}(x)$ for $i=1,2$ the deterministic flow defined to be the solution of the control problem. 
$\psi^{(i)}_{st}(x)=x+\int_{s}^{t}V_{i}\left(\psi^{(i)}_{sr}(x)\right)dr$.
\begin{proposition}\label{be:nospeed}
For $t\leq t_{u}^{n}$, $z\in\tilde{U}^{n}_{Q}$ we have 
$\left\|\psi^{(1/2)}_{0,t}(z)-z-tw_{1/2}\right\|\leq \epsilon t$.
\end{proposition}
Proof of Proposition~\ref{be:nospeed}: For $z\in \tilde{U}^n_Q$ we have $\left\|V_{1/2}(z)\right\|\leq 0,5n\epsilon \left(t_{u}^{n}\right)^{-1}$, which implies for $z\in\tilde{U}_{Q}^{n}$ that
\be  
\inf_{t>0}\left\{\left\|\psi_{0t}^{(i)}(z)-z\right\|>\frac{n\epsilon}{2}\right\}\geq t_{u}^{n}\Rightarrow\nonumber\\
\inf_{t>0}\left\{\psi_{0t}^{(i)}(z)\notin U_{Q}^{n} \right\}\geq t_{u}^{n},
\ee  
which proves Proposition~\ref{be:nospeed} because we have for $z\in \tilde{U}^n_Q$ that
\bee  
\left\|\psi_{0t}^{(i)}(z)-z-tw_{i}\right\|\leq\int_{0}^{t}\left\|V_{i}\left(\psi_{0s}^{(i)}(z)\right)-w_{i}\right\|ds\leq\epsilon t.
\eee  
\hfill$\Box$\\
Considering the following coordinates: $Z:=(Z_{1},Z_{2}):\tilde{U}_{Q}^{102}\rightarrow]-51,51[$, generated by the  vectorfields $\left(\epsilon w_{i}:i=1,2\right)$ we choose $U_{Q}:=Z^{-1}\left(]-1,1[^{2}\right)$.
\subsubsection{From Large To Positive Probability}
As we will see it suffices to show that there is $0<\theta<1$ such that for a large curve wtih a non-empty intersection with $K_{R}(P)$ we have uniformly in  $Q\in K_{R}(P)$ that
\begin{equation}\label{eq:cp}
\cp{\tilde{\tau}_{Q}<t_{j}}{\tilde{\tau}_{Q}>t_{j-1}}\geq \theta.
\end{equation}
Therein for a $T>0$ (to be specified later) let $t_{0}:=0$ and for $j\in\mathbb{N}$ define \[t_{j}:=\inf\left\{t\in\mathbb{R}:t\geq t_{j-1}+1+T:\gamma_{t}\cap K_{R}(P)\neq\emptyset,\diam(\gamma_{t})\geq1 \right\}.\]
Following Theorem~\ref{th:cor6} there is $C_{14}>0$ and for $m\in\mathbb{N}$ a $\kappa_m^{(10)}\in\mathbb{R}$, such that for $j\in\mathbb{N}$ we have $\prob{t_{j}>C_{14}j}\leq \kappa_m^{(10)}j^{-m}$  which implies 
\begin{align}
&\prob{\tilde{\tau}_{Q}>t}
\leq\prob{t_{\left\lfloor\frac{t}{C_{14}}\right\rfloor}>t}+\prob{\tilde{\tau}_{Q}> t_{\left\lfloor\frac{t}{C_{14}}\right\rfloor}}\nonumber\\
\leq&\prob{t_{\left\lfloor\frac{t}{C_{14}}\right\rfloor}>C_{14} \left\lfloor\frac{t}{C_{14}}\right\rfloor}+\cp{\tilde{\tau}_{Q}> t_{\left\lfloor\frac{t}{C_{14}}\right\rfloor}}{\tilde{\tau}_{Q}>t_{\left\lfloor\frac{t}{C_{14}}\right\rfloor-1}}\prob{\tilde{\tau}_{Q}>t_{\left\lfloor\frac{t}{C_{14}}\right\rfloor-1}}\nonumber\\
\leq&\prob{t_{\left\lfloor\frac{t}{C_{14}}\right\rfloor}>C_{14} \left\lfloor\frac{t}{C_{14}}\right\rfloor}+\ldots
\leq\kappa_m^{(10)}\left\lfloor\frac{t}{C_{14}}\right\rfloor^{-m}+(1-\theta)^{\left\lfloor\frac{t}{C_{14}}\right\rfloor}
\leq \kappa_m^{(9)}t^{-m}\nonumber
\end{align}
for suitable $\kappa_m^{(9)}\in\mathbb{R}$. So we have only to prove (\ref{eq:cp}). Therefore it is enough to show that for $\gamma$ (as before) there are $T>0$ and $\theta>0$ (not depending on the chosen $\gamma$) such that we have uniformly in $Q\in K_{R}(P)$ that
\begin{equation}\label{eq:prob}
\prob{U_{Q}\subset\bigcup_{0\leq s\leq T}\gamma_{s}}\geq\theta.
\end{equation}
\subsubsection{Approaching The Small Square}
\begin{figure}\centering \resizebox{3cm}{!}{\includegraphics{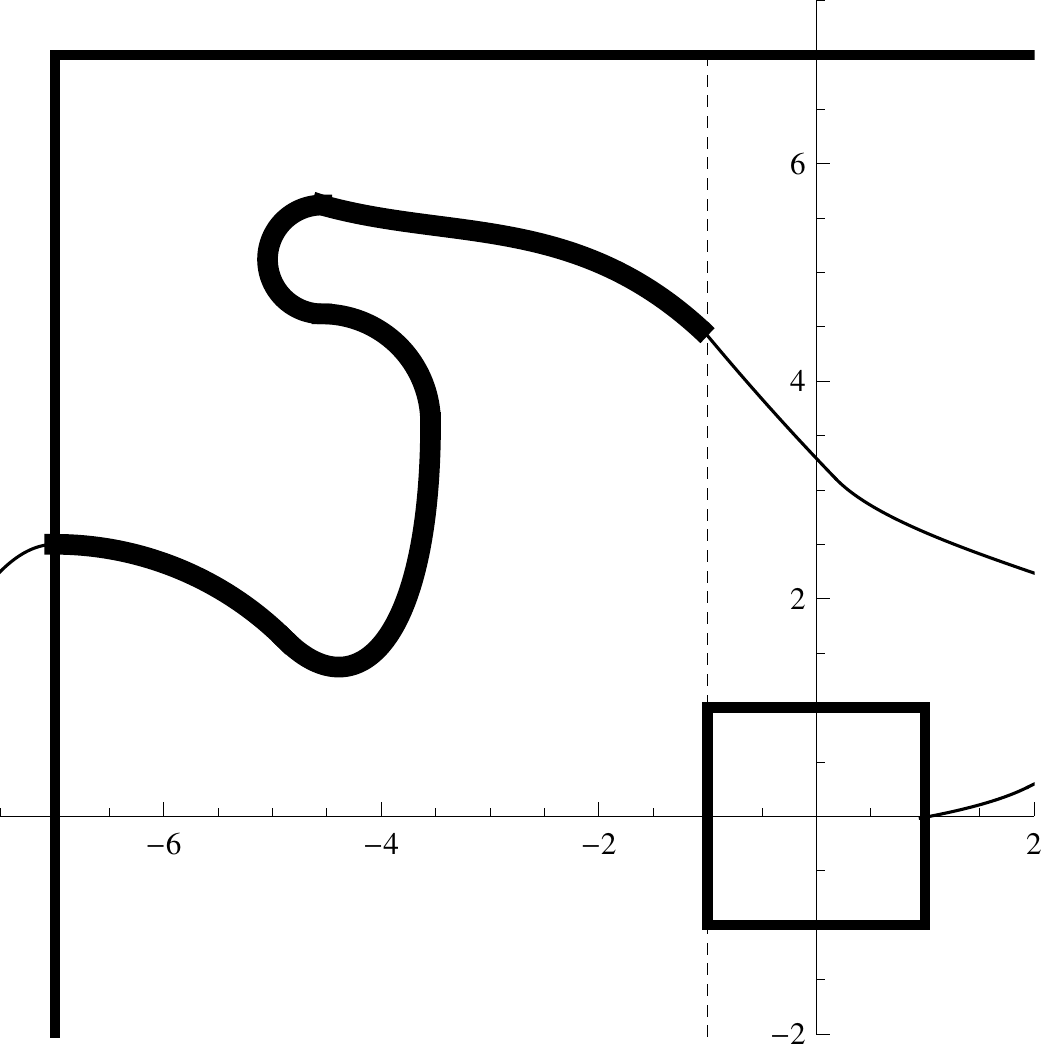}}\caption{$\tilde{\gamma}$ (fat), rest of $\hat{\gamma}$ (regular), The endpoints of $\tilde{\gamma}$ are $z$ (left) and $\tilde{y}$ (right). \label{abb:vorsweep}}\end{figure}

Let $\hat{U}_{Q}:=Z^{-1}\left(]-7,7[^{2}\right)$. Then we have obviously $U_{Q}\subset\hat{U}_{Q}$. Choose $x$ and $y$ in $\gamma$ with $\left\|x-P\right\|\leq R$ and $\left\|x-y\right\|\geq 0.5$.
Due to the Lemmas~\ref{le:eigenwerte} and \ref{le:korrfunk} the eigenvalues of $\bar{b}^{2}(z):=\bar{b}^{*}(z)\bar{b}(z)$  are bounded below by a positive constant $C_{4}$ on $\left\{\left\|z\right\|\geq\delta\right\}$ for arbitrary $\delta>0$. The boundedness of the correlation functions gives an upper bound $C_{5}$. Therefore the $\mathbb{R}^{4}$-valued semimartingale 
\[\left\{\left(\begin{array}{l}x_{t}-\left(x+2t(Q-x)\right)\\y_{t}-\left(y+2t(Q-x)\right) \end{array}\right):t\in\left[0,\frac{1}{2}\right]\right\}\] satisfies the assumptions of~\cite[Lemma 2.4]{ss}. So this lemma yields for $t=0.5$ and $\delta=0.5\epsilon$  ($C_{4}$, $C_{5}$ and $C_6$ can be chosen to be independent of $x$ and $y$): 
\begin{equation}\prob{x_{\frac{1}{2}}\in U_{Q},y_{\frac{1}{2}}\notin \hat{U}_{Q}}\geq p,
\end{equation} where $p>0$ does not depend on the special choice of $\gamma$, because $\epsilon\leq(56\sqrt{2})^{-1}$ implies $\diam\hat{U}_{Q}\leq0.25$.
Denote by $\hat{\hat{\gamma}}$ the subcurve of $\gamma$, between $x_{0.5}$ and $y_{0.5}$ and by $\hat{\gamma}$ a minimal subcurve of $\hat{\hat{\gamma}}$, which is contained in $\hat{U}_{Q}$ and which links $\partial\hat{U}_{Q}$ to $\partial U_{Q}$ (minimal means that no proper subcurve has these properties). Due to minimality of $\hat{\gamma}$ the set $\hat{\gamma}\cap \partial \hat{U}_{Q}$ consists of a single point which we will denote by $z$.
$\partial\hat{U}_{Q}$ consists of four pieces. Without loss of generality assume $z\in Z^{-1}\left(\{-7\}\times [-7,7]\right)$ (the other cases are similar).Let $\tilde{\gamma}$ be the minimal subcurve of $\hat{\gamma}$, linking $z$ with $Z^{-1}(\{-1\}\times[-7,7])$ and $\tilde{y}:=\tilde{\gamma}_{0.5}\cap Z^{-1}(\{-1\}\times[-7,7])$ the intersection point (see~fig.~\ref{abb:vorsweep}). 
We have to show that there exist $T>0$ and $\theta>0$ such that for any curve $\tilde{\gamma}\subset Z^{-1}\left([-7,7]\times[-7,7]\right)$ linking $Z^{-1}(\{-7\}\times[-7,7])$ to $Z^{-1}(\{-1\}\times[-7,7])$ the following holds: \begin{equation}\label{eq:tilga}\prob{U_{Q}\subset\cup_{0\leq s\leq T}\tilde{\gamma}_{s}}\geq\theta.\end{equation}
\subsubsection{Reduction To A Control Problem}
For an $\mathcal{H}$-simple control $V$ denote by $\psi_{s}^{(V)}(x)$ the solution to the control problem
$$\left\{
\begin{array}{ccc}
\partial_{t}\psi_{t}^{(V)}(x)&=&V\left(\psi^{(V)}_{t}(x)\right)\\
\psi_{0}^{(V)}(x)&=&x
\end{array}\right..$$
Assume we can construct a $\mathcal{H}$-simple control $V$ with the following property:
If $\Psi(.,.)$ is a continuous mapping from $\left[0,T\right]\times\mathbb{R}^{2}$ to $\mathbb{R}^{2}$  which satisfies
\begin{equation}\label{eq:sweepapp}
\left|Z\left(\Psi(s,x)\right)-Z\left(\psi^{(V)}_{s}(x)\right)\right|<\frac{1}{2}\end{equation} for $x\in\tilde{\gamma}$ and $s\in[0,T]$, then we also have
$U_{Q}\subset\bigcup_{x\in\tilde{\gamma}}\bigcup_{0\leq s\leq1}\Psi(s,x)$.
Then Theorem~\ref{sa:ibfsupp} applied to the intervals of constance of $V$ and the independence properties of Brownian flows prove (\ref{eq:tilga}).
Let $\tilde{\gamma}=\left\{\tilde{\gamma}(u):u\in[0,1]\right\}$ be a parametrization of $\tilde{\gamma}$ and $\tilde{\psi}(s,u):=\psi^{(V)}_{s}\left(\tilde{\gamma}(u)\right)$ as well as  $\tilde{\Psi}(s,u):=\Psi(s,\tilde{\gamma}(u))$.\\
($\bar{\psi}$ and $\bar{\Psi}$ are defined similarly with $\tilde{\gamma}$ replaced by $\bar{\gamma}$, see its definition below.)
We want to construct a $\mathcal{H}$-simple control $V$, that implies for any $\Psi$ fulfilling (\ref{eq:sweepapp}) that $U_{Q}\subset \cup_{\frac{7}{17}T\leq s\leq T}\cup_{0\leq u\leq1}\tilde{\Psi}(s,u)$.\\
Set $\bar{\gamma}:=\partial\left(\left[\frac{7}{17}T,T\right]\times\left[0,1\right]\right)$. $V$ is supposed to yield for $\tilde{Q}\in U_{Q}$ that \begin{equation}\label{eq:index}\textnormal{ind}\left(\bar{\Psi},\tilde{Q}\right)=1.\end{equation}
Therein we denote by $\textnormal{ind}\left(\bar{\Psi},\tilde{Q}\right)$ the curving number of $\bar{\Psi}$ arround $\tilde{Q}$. To show (\ref{eq:index}) for all $\tilde{\Psi}$ with \begin{equation}\label{eq:indexvor}\left\|Z(\tilde{\Psi}(.,.))-Z(\tilde{\psi}(.,.))\right\|_{\infty}\leq 0.5\end{equation} we construct $V$ in a way, that provides for $\tilde{Q}\in U_{Q}$ the following
\begin{equation}\label{eq:index2}
\textnormal{ind}\left(\bar{\psi},\tilde{Q}\right)=1 \textnormal{ and } \dist\left(Z\left(\bar{\psi}\right),Z(U_{Q})\right)\geq 1.
\end{equation}
Note that (\ref{eq:index2}) implies that for any $\tilde{\Psi}$ satisfiing (\ref{eq:indexvor}) we indeed get (\ref{eq:index}). ($\tilde{\psi}$ sweeps the entire set $Z^{-1}\left([-2,2]^{2}\right)$).
\subsubsection{Construction Of A Sweeping Control}
\begin{figure}\centering \resizebox{4cm}{!}{\includegraphics{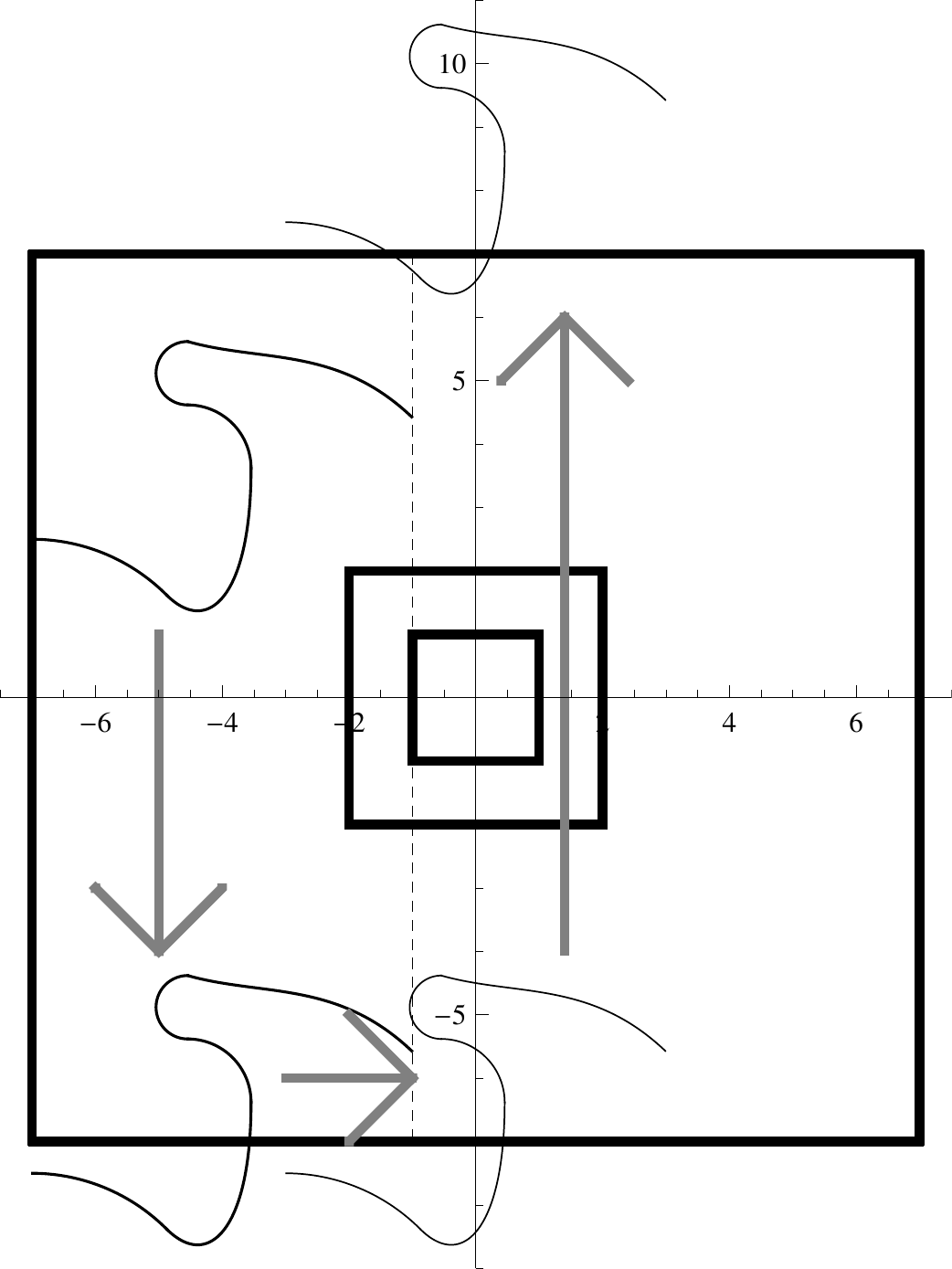}}\caption{$\gamma$ (curve), $\partial U_{Q}$ and $Z^{-1}\left(\{-1\}\times[-7,7]\right)$ (dashed) , $Z^{-1}([-2,2]^{2})$ and $\partial\hat{U}_{Q}$ (black), way of $\gamma$ along the control (gray) \label{abb:sweep}}\end{figure}

Consider the following $\mathcal{H}$-simple control $V:[0,T]:=[0,34\epsilon]\rightarrow\mathcal{H}$:
\begin{equation}
V(.,t):=\left\{
\begin{array}{rl}
-V_{2}(.)&:t\in[0,10\epsilon[\\
V_{1}(.)&:t\in[10\epsilon,14\epsilon[\\
V_{2}(.)&:t\in[14\epsilon,34\epsilon]
\end{array}\right..
\end{equation}
This control satisfies all our wishes (see~Fig.~\ref{abb:sweep}).
Note that $34\epsilon=\frac{102}{3}\epsilon\leq t^{102}_{u}$ ensures the suitability of Proposition~\ref{be:nospeed} for $t\leq34\epsilon$. So we get for $t\leq34\epsilon$, $z\in\hat{U}_{Q}$ and $i=1,2$ that
$\left\|\Psi^{(i)}_{0t}(z)-z-tw_{i}\right\|\leq 34\epsilon^{2}\leq\frac{\epsilon}{3}$,
which gives the following estimates for $x\in\tilde{\gamma}$.
\begin{align}  
-7&\leq& Z_{1}(x)&\leq& -1;&&&-7&\leq& Z_{2}(x)&\leq& 7&:&t=0\nonumber,\\
-7-\frac{1}{3}&\leq& Z_{1}(x_{t})&\leq& -1+\frac{1}{3};&&&-17-\frac{1}{3}&\leq& Z_{2}(x_{t})&\leq& -3+\frac{1}{3}&:&t=10\epsilon\nonumber,\\
-3-\frac{2}{3}&\leq& Z_{1}(x_{t})&\leq& 3+\frac{2}{3};&&&-17-\frac{2}{3}&\leq& Z_{2}(x_{t})&\leq& -3+\frac{2}{3}&:&t=14\epsilon\nonumber,\\
-3-1&\leq& Z_{1}(x_{t})&\leq& 3+1;&&&3-1&\leq& Z_{2}(x_{t})&\leq& 17+1&:&t=34\epsilon\nonumber.
\end{align}
Herein we used that elapsing of time through the intervals of constance of $V$ changes any coordinate as if $Z$ was generated by $\epsilon V$ with an error of at most $3^{-1}$ (see Fig.~\ref{abb:sweep}). A similar arguing holds for $z$ and $\tilde{y}$.
\begin{align}
&&Z_{1}(z)&=&-7;&&&&Z_{1}(\tilde{y})&=&-1&:&t=0\nonumber,\\
-7-\frac{1}{3}&\leq& Z_{1}(z_{t})&\leq& -7+\frac{1}{3};&&-1-\frac{1}{3}&\leq& Z_{1}(\tilde{y}_{t})&\leq& -1+\frac{1}{3}&:&t=10\epsilon\nonumber,\\
-3-\frac{2}{3}&\leq& Z_{1}(z_{t})&\leq& -3+\frac{2}{3};&&3-\frac{2}{3}&\leq& Z_{1}(\tilde{y}_{t})&\leq& 3+\frac{2}{3}&:&t=14\epsilon\nonumber,\\
-3-1&\leq& Z_{1}(z_{t})&\leq& -3+1;&&3-1&\leq& Z_{1}(\tilde{y}_{t})&\leq& 3+1&:&t=34\epsilon\nonumber.
\end{align}  
This means that for $t\in\left[14\epsilon,34\epsilon\right]$ $z_{t}$ is on the left and $\tilde{y}_{t}$ is on the right of $\check{U}$  which implies (\ref{eq:index2}) and completes the proof of Lemma~\ref{le:sweep}. \hfill$\Box$\\ 
\subsection{Dependence Of $\left\|v\right\|^{R}$ On $R$}
>From now on we leave the ideas of~\cite{dkk} and show the following directly:
\begin{lemma}\label{le:const}
If we define for $R>0$, $\tilde{R}\geq1$ and $v\in\mathbb{R}^d$ $\left\|v\right\|_{\tilde{R}}^{R}$ via 
$\left\|v\right\|^{R}_{\tilde{R}}:=\lim_{t\rightarrow\infty}\frac{ \sup_{\gamma\in\mathcal{C}_{\tilde{R}}}\expec{\tau^{R}(\gamma,tv)} }{t}$, 
then this limit exists and we have for arbitrary $R_{1}>0$, $R_{2}>0$ and $\tilde{R}_{1}\geq1$, $\tilde{R}_{2}\geq1$ that $\left\|.\right\|^{R_{1}}_{\tilde{R}_{1}}\equiv\left\|.\right\|^{R_{2}}_{\tilde{R}_{2}}$ i.e. $\left\|v\right\|^{R}$ does not depend on $R$.
\end{lemma}
Proof: Define for $t\geq 0$, $R>0$ and $\tilde{R}\geq1$ the function $\bar{g}=\bar{g}(R,\tilde{R},t)$ via $\bar{g}(R,\tilde{R},t):=\sup_{\gamma\in\mathcal{C}_{\tilde{R}}}\left\{\expec{\tau^{R}(\gamma, tv)} \right\}$.
Herein fix $v\in\mathbb{R}^d$ with $\left\|v\right\|=1$. As already seen there is $R>0$ such that 
\begin{equation}\label{eq:rrkonv} \left\|v\right\|_{R}^{R}:=\lim_{t\rightarrow\infty}\frac{g(R,R,t)}{t} \end{equation}
exists (the limit was named $\left\|v\right\|^{R}$). 
We will first prove that if
we fix $R\geq1$ in a way that we have convergence in (\ref{eq:rrkonv}), then we have for arbitrary $\tilde{R}\geq1$ that
\begin{equation}\label{eq:stattprop1}
\left\|v\right\|_{\tilde{R}}^{R}=\lim_{t\rightarrow\infty}\frac{g(R,\tilde{R},t)}{t}=\left\|v\right\|_{R}^{R}.
\end{equation}
\newpage
Observe:
\begin{enumerate}
\item If $\tilde{R}\geq R$ then $\bar{g}(R,\tilde{R},t)\geq \bar{g}(R,R,t)$ is obvious.\\ Isotropy yields:
$\bar{g}(R,\tilde{R},t)\leq\bar{g}(R,R,t+\tilde{R})$.
For we have that\\ $\expec{\tau^{R}\left(\gamma,(t+\tilde{R})v\right)}\leq\expec{\tau^{R}(\gamma,tv)}+C_{15}$
for a $C_{15}>0$ (uniformly chosen in $\gamma$) we get
$\bar{g}(R,R,t)\leq \bar{g}(R,\tilde{R},t)\leq \bar{g}(R,R,t)+C_{15}.$
\item If $\tilde{R}<R$ we obtain similarly that $\bar{g}(R,\tilde{R},t)\leq \bar{g}(R,R,t)\leq\bar{g}(R,\tilde{R},t)+C_{15}$
\end{enumerate}
Sending $t\rightarrow\infty$ proves~\eqref{eq:stattprop1} from the latter.\\
Now we will prove that
$\left\|v\right\|^{R}_{1}$ exists for any $R>0$ and we have:
\begin{equation}\label{eq:stattprop2}
\left\|.\right\|_{1}^{R}\equiv\left\|.\right\|_{1}^{\tilde{R}}
\end{equation}
Without loss of generality assume that $\tilde{R}>R$. Then $\bar{g}(\tilde{R},1,t)\leq \bar{g}(R,1,t)$ is obvious. Addionally (in $\tau^{R}(\gamma,tv)$ one takes a subcurve if necessary) we have 
$\expec{\tau^{R}(\gamma,tv)}\leq\expec{\tau^{\tilde{R}}(\gamma,tv)}+\sup_{\check{\gamma}\in\mathcal{C}_{1}}\expec{\tau^{R}(\check{\gamma},\tilde{R}v)}$,
so~\eqref{eq:stattprop2}  follows via $t\rightarrow\infty$ from
$\bar{g}(R,1,t)\leq\bar{g}(\tilde{R},1,t)+C_{16}$
for some $C_{16}>0$.
The proof of Lemma~\ref{le:const} is complete.\hfill$\Box$\\
\section{The Upper Bound}
\label{sec:ub}
\subsection{The Speed Of A Slow Curve Asymptotically Has A Dirac Distribution On The Proper Time Scale }
Let $R$ be as before and $v\in\R^d$ with $\left\|v\right\|^R=1$. Choose for any $t\in\left]0,\infty\right[$ a curve $\gamma^{(t)}\in\mathcal{C}_{R/2}$ such that $|\expec{\tau^R(\gamma^{(t)},tv)}-|tv|^R|t^{-1}\to0$.
By the definition of $\left\|v\right\|^R$ we already know
\be\label{eq:expecconv}
\left|\frac{ \expec{\tau^R(\gamma^{(t)},tv)} }{t}-\left\|v\right\|^R \right|\leq \left|\frac{\expec{\tau^R(\gamma^{(t)},tv)}-|tv|^R}{t}\right|+\left|\frac{|tv|^R}{t}-\left\|v\right\|^R \right| \rightarrow0. 
\ee
We will investigate the asymptotic law of $\tau^R(\gamma^{(t)},tv)$ in the following lemma.
\begin{lemma}\label{le:diracdist}
Let $\left(X_t:t>0 \right)$ be a family of integrable random variables such that we have 
$\lim_{t\rightarrow\infty}\expec{X_t-\expec{X_t};X_t-\expec{X_t}>\delta}=0$ for any $\delta>0$.
Then we also have for $\delta>0$ that $\lim_{t\rightarrow\infty}\prob{X_t-\expec{X_t}<-\delta}=0$. 
\end{lemma} 
Proof: Denote $X_t-\expec{X_t}$ by $Y_t$ and suppose that we can find $\delta>0$, $\epsilon>0$ und a sequence $(t_n)_{n\in\mathbb{N}}$ such that for $t_n\nearrow\infty$ and any $n$ we have $\prob{Y_{t_n}<-\bar{\delta}}>\epsilon$. This immediately yields\\
$0=\expec{Y_{t_n}}=\expec{Y_{t_n};Y_{t_n}\leq-\delta}
+\expec{Y_{t_n};Y_{t_n}\in(-\delta,\frac{\epsilon\delta}{2})}+\expec{Y_{t_n};Y_{t_n}\geq\frac{\epsilon\delta}{2}} $ i.e. a contradiction for large $n$ because the right-hand site of the latter is strictly negative for large $n$.\hfill$\Box$\\
\begin{lemma}
Let $R$ be as before and fix $v$ with $\left\|v\right\|^R=1$, Then for any $\delta>0, n,m\in\mathbb{N}$ there are $K_m^{(n)}(\delta)\in\R$ and $\tilde{t}_m^{(n)}(\delta)$ such that we have for all $t\geq0$ that
\be
\expec{\left(\tau^R(\gamma^{(t)},tv)t^{-1}\right)^n;\tau^R(\gamma^{(t)},tv)t^{-1}>(1+\delta)}\leq K_m^{(n)}(\delta)t^{-m}\nonumber
\ee
and for $t\geq\tilde{t}_m^{(n)}(\delta)$ that
\be
\expec{\left(\tau^R(\gamma^{(t)},tv)t^{-1}\right)^n;\tau^R(\gamma^{(t)},tv)t^{-1}>\expec{\tau^R(\gamma^{(t)},tv)t^{-1}}+\delta}\leq K_m^{(n)}(\delta)t^{-m}.\nonumber
\ee
\end{lemma}
Proof: The first equation follows from straightfoward estimates using~\eqref{eq:hittail} and the second one is implied by the first one and \eqref{eq:expecconv}.\hfill$\Box$\\
The previous lemma implies that we can apply Lemma~\ref{le:diracdist} to\\ $X_t:=\tau^R(\gamma^{(t)},tv)t^{-1}$ to conclude that $X_t$ converges to $1$ in probability.
\subsection{Time Reverse - Comparison Of Fast And Slow Curves}
\begin{figure}\centering  \resizebox{6cm}{!}{\includegraphics{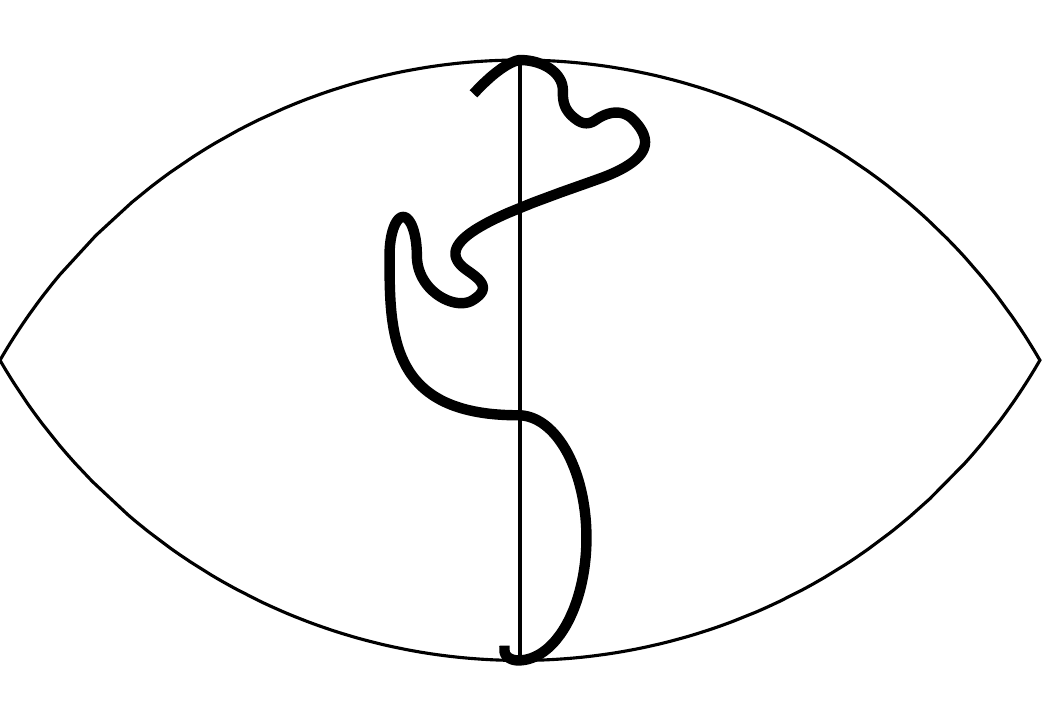}}\resizebox{4cm}{!}{\includegraphics{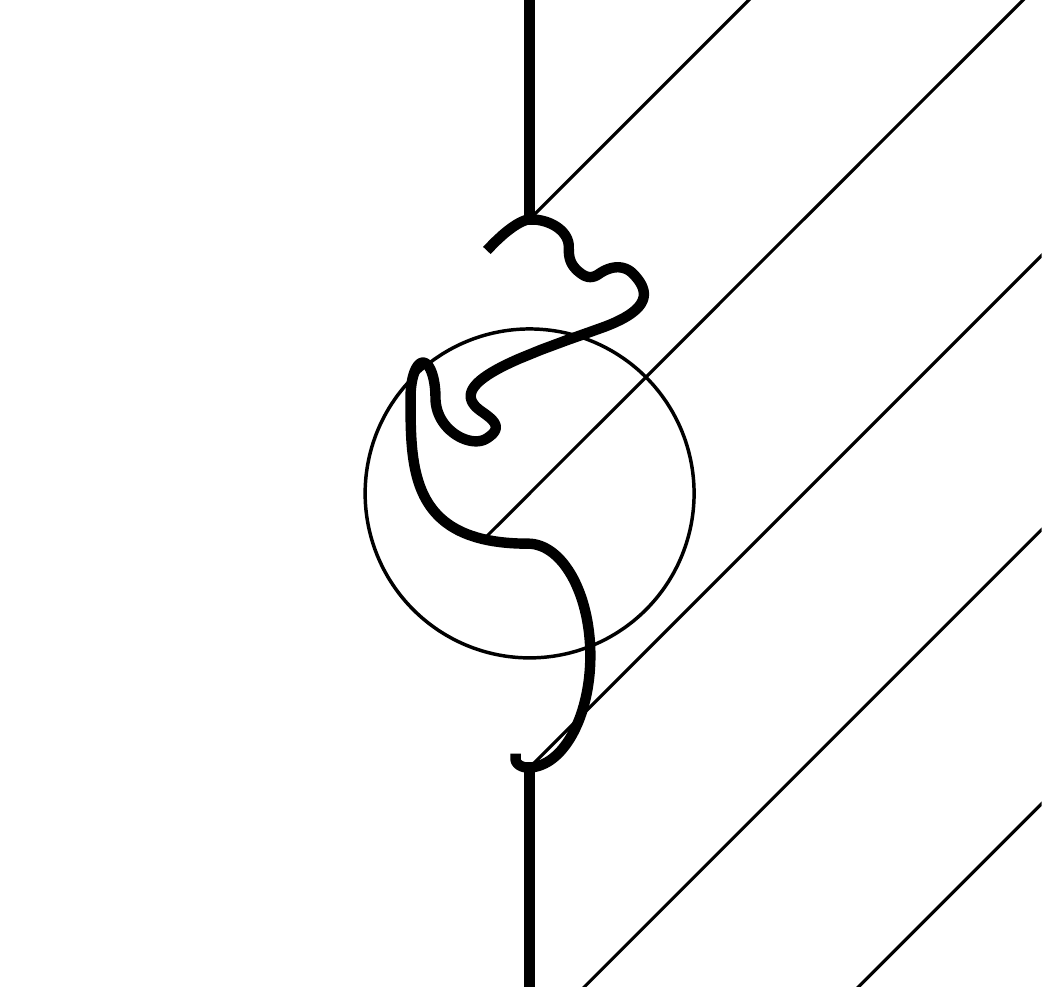}}
\caption{Definition of $\gamma_{t_1}$ via intersection of circles arround its endpoints  \label{abb:kurve}}
\caption{ Adding to half lines to $\gamma_{t_1}$ divides $\R^2$ into black and white parts \label{abb:teilung}}\end{figure}

We will have to assume $d=2$ from now on (unless otherwise stated) for the following arguments strongly depend on the topology of the plane. 
\begin{theorem}\label{th:allesdirac}
Let $\Gamma:=\partial K_R(0)$. There are $C_{20}>0$ and $C_{17}>0$ and for $m\in\mathbb{N}$ $\kappa_m^{(11)}\in\R$ such that we have for $T\geq\sqrt{t}$ and any $\gamma\in\mathcal{C}_{R/2}$ that
\be
\kappa_m^{(11)}T^{-m}+\prob{\tau^{R}(\gamma,tv)\leq T+C_{17}}\geq C_{20}\prob{\tau^{R}(\Gamma,tv)\leq T}.
\ee
\end{theorem}
The proof of Theorem~\ref{th:allesdirac} uses the following lemmas. Denote by $\mathcal{C}_R^*$ the set of all large curves $\gamma$ with $\gamma\cap\partial K_R(0)\neq\emptyset$.
\begin{lemma}\label{le:q1}
There is a constant $C_{18}>0$ with\\
$
\inf_{\gamma\in\mathcal{C}_{R}^*}\inf_{t\geq C_{18}}\prob{\gamma_t\cap\partial K_{R}(0)\neq\emptyset;\diam(\gamma_t)\geq1}=:p_1>0.
$
\end{lemma}
Proof: Since for any constant $C_{18}$ the distance process $\dist(\gamma_{t-0.5C_{18}},0)$ of a long curve from the origin can be majorized by a stationary process there is $\epsilon>0$ with
$\inf_{\gamma\in\mathcal{C}_{R}^*}\inf_{t\geq C_{18}}\prob{\gamma_t\cap K_{R}(0)\neq\emptyset}\geq\epsilon$. 
Now choose $C_{18}$ large enough for
$\inf_{\gamma\in\mathcal{C}_{R}^*}\inf_{t\geq C_{18}}\prob{\textnormal{diam}\gamma_t>3R}\geq1-\frac{\epsilon}{2}$ to hold. This directly implies
\[
\inf_{\gamma\in\mathcal{C}_{R}^*}\inf_{t\geq C_{18}}\prob{\gamma_t\cap K_R(0)\neq\emptyset;\gamma_t\cap K_R(0)^C\neq\emptyset;\diam(\gamma_t)\geq1}\geq\frac{\epsilon}{2}
\] 
and hence we get $p_1>0$ completing the proof of Lemma~\ref{le:q1}. \hfill$\Box$\\
\begin{lemma}\label{le:q2}
There is $C_{19}>0$ with 
$\inf_{\bar{\gamma}\in\mathcal{C}_R^*}\inf_{\gamma\in\mathcal{C}_{R/2}}\prob{\bar{\gamma}_{C_{19}}\cap\gamma\neq\emptyset}=:p_2>0$.
\end{lemma}
Proof: First write $t=t_{1}+t_{2}$ for some non-negative $t_{1}$, $t_{2}$ and observe
\begin{align}
&\prob{\gamma_{t}\cap\bar{\gamma}=\emptyset}&=&\prob{\Phi_{t_{1},t}\left(\Phi_{0,t_{1}}(\gamma)\right)\cap\bar{\gamma}=\emptyset}\nonumber\\
=&\prob{\Phi_{0,t_{1}}\left(\gamma\right)\cap\Phi_{t,t_{1}}\left(\bar{\gamma}\right)=\emptyset}&=&\left(\mathbb{P}\otimes\mathbb{P}\right)\circ\pi_{1}^{-1}\left[\Phi_{0,t_{1}}\left(\gamma\right)\cap\Phi_{t,t_{1}}\left(\bar{\gamma}\right)=\emptyset\right]\nonumber\\
=&\mathbb{P}\otimes\mathbb{P}\left[\Phi_{0,t_{1}}\left(\gamma\right)\cap\tilde{\Phi}_{t,t_{1}}\left(\bar{\gamma}\right)=\emptyset\right]&=&\mathbb{P}\otimes\mathbb{P}\left[\Phi_{0,t_{1}}\left(\gamma\right)\cap\tilde{\Phi}_{t_{1},t}\left(\bar{\gamma}\right)=\emptyset\right]\nonumber\\
=&\mathbb{P}\otimes\mathbb{P}\left[\Phi_{0,t_{1}}\left(\gamma\right)\cap\tilde{\Phi}_{0,t_{2}}\left(\bar{\gamma}\right)=\emptyset\right].
\end{align}
Herein we denote by $\tilde{\Phi}$ an independent copy of $\Phi$ (for example defined on  $\left(\Omega\times\Omega,\mathcal{F}\otimes\mathcal{F},\mathbb{P}\otimes\mathbb{P} \right)$). So we can (instead of having $\bar{\gamma}$ running $t$) split $t$ and let $\gamma$ run one part of the time and $\bar{\gamma}$ the rest of it. We already know that for sufficiently large $t_1$ $\gamma_{t_1}$ is dense as much as $\sqrt{t_1}$ in $K_{t_{1}^{0.9}}(0)$ with probability say at least $0.5$ uniformly in $\gamma$ (this probability converges to one as $t_1\rightarrow\infty$). Rename $\gamma_{t_1}$ to be a connected subcurve $\gamma_{t_1}\cap K_{t_{1}^{0.9}}(0)$ of diameter $t_{1}^{0.8}$ that has distance not more than $\sqrt{t_1}$ from the origin. We may asume that the endpoints of the new $\gamma_{t_1}$ have distance $t_{1}^{0.8}$ from each other (which we do). Note that $\gamma_{t_1}$ is contained in the intersection of the $t^{0.8}_1$-balls around its endpoints.(see Fig.~\ref{abb:kurve})
By adding two half lines to $\gamma_{t_1}$ we cut the plane into two parts (see Fig.~\ref{abb:teilung}) 
say the black part and the white part.  Now we fix a ball $\mathcal{K}$ of radius of $t_1^{0.6}$ centered on the perdendicular bisector of the endpoints of $\gamma_{t_1}$ exactly one half of which is black (measured with Lebesgue measure). The existence of such a ball follows from a continuity argument and the fact that there are completely black  and completely white balls centered there. Fix $t_1$ large enough for $t^{0.9}_1>>t^{0.8}_1>>t^{0.7}_1>>t^{0.6}_1>>t^{0.5}_1>>1$ to hold. Observe now that any curve that links the black part to the white part of $\mathcal{K}$ without intersecting $\gamma_{t_1}$ must have diameter at least $t^{0.7}_1$. Of course all the choices above can be made $\mathcal{F}_{t_1}$-mesureable. \\
Now it is $\bar{\gamma}$'s turn to do the rest within time $t_2=t_3+1$. Choose a point in $\gamma_{t_1}$ that has distance at least $t^{0.5}$ to the complement of $\mathcal{K}$ such that at least one third of its $2R$-neighbourhood is black and white respectively.  Fix $t_3$ large enough that the probability of the event that  $\tilde{\Phi}(\bar{\gamma}_{t_3})$ has distance to this point less or equal to $R$ (and $\tilde{\Phi}(\bar{\gamma}_{t_3})$ is long) is at least $\epsilon$ uniformly in $\bar{\gamma}$. 
So with probability at least $0.5\epsilon$ a point say $x$ in $\tilde{\Phi}(\bar{\gamma}_{t_3})$  has an environment of diameter $3R$ at least  one percent of which is white and black respectively. Choose another point in $\tilde{\Phi}(\bar{\gamma}_{t_3})$ say $y$ with distance $1/2$ of $x$ such that the subcurve (denoted $\check{\gamma}$) of $\tilde{\Phi}(\bar{\gamma}_{t_3})$ linking $x$ and $y$ has diameter $1/2$ and observe now that the lemma follows from Theorem~\ref{th:posdens} because we can choose $t_1$ large enough for $\check{\gamma}$ not to reach $\mathcal{K}^C$ within the remaining time $1$ with sufficiently large probability. With $C_{19}:=t_1+t_3+1$ the proof is complete.\hfill$\Box$\\
\newpage 
\begin{lemma}\label{le:q3}
There is $p_3>0 $ such that we have\\
$\prob{\tau^R(\Gamma,tv)\leq T}\leq \frac{1}{p_3}\prob{\tau^R(\Gamma,tv)\leq T;\diam(\Gamma_{T+C_{18}})\geq 1}$.
\end{lemma}
Proof: This a direct consequence of the fact, that the diameter of long curves uniformly has a chance to grow to infinity without being smaller than $1$ after time $C_{18}$. If necessary we increase $C_{18}$ (without changing notation).\hfill$\Box$\\  
Proof of Theorem~\ref{th:allesdirac}: First we use Lemma~\ref{le:q1} to estimate
\bee 
&&\prob{\Gamma_{T+C_{18}}\cap\partial K_R(tv)\neq\emptyset}\geq\prob{\tau^R(\Gamma,tv)\leq T;\diam(\Gamma_{T+C_{18}})\geq1}\nonumber\\
&\cdot&\cp {\diam(\Gamma_{T+C_{18}})\geq 1;\Gamma_{T+C_{18}}\cap\partial K_R(tv)\neq\emptyset}{\tau^R(\Gamma,tv)\leq T}\nonumber
\eee 
which implies with Lemma~\ref{le:q3}
\begin{align}
&\prob{\tau^{R}(\Gamma,tv)\leq T}\leq \frac{1}{p_3}\prob{\tau^R(\Gamma,tv)\leq T;\diam(\Gamma_{T+C_{18}})\geq 1}\nonumber\\
\leq&\frac{1}{p_1p_3} \prob{\Gamma_{T+C_{18}}\cap\partial K_R(tv)\neq\emptyset}\nonumber\\
\leq&\frac{1}{p_1p_3}\left( \prob{\Gamma_{T+C_{18}}\cap\partial K_R(tv)\neq\emptyset;\diam(\Gamma_{T+C_{18}})\geq1}+\prob{\diam\Gamma_{T+C_{18}}<1 }\right).\nonumber
\end{align}
Usage of Lemma~\ref{le:reverse} and symmetry that $\prob{\tau^{R}(\Gamma,tv)\leq T}p_1p_3 $ is at most
\begin{align}
&\prob{\Gamma\cap\partial K_R(tv)_{T+C_{18}}\neq\emptyset;\diam(\partial K_R(tv)_{T+C_{18}})\geq1}\nonumber\\
\leq&\frac{\prob{\gamma\cap\partial K_R(tv)_{T+C_{18}+C_{19}}\neq\emptyset;\Gamma\cap\partial K_R(tv)_{T+C_{18}}\neq\emptyset}}{\cp{\gamma\cap\partial K_R(tv)_{T+C_{18}+C_{19}}\neq\emptyset}{\Gamma\cap\partial K_R(tv)_{T+C_{18}}\neq\emptyset;\diam(\partial K_R(tv)_{T+C_{18}})\geq1}}\nonumber\\
\leq&\frac{\prob{\gamma\cap\partial K_R(tv)_{T+C_{18}+C_{19}}\neq\emptyset}}{\cp{\gamma\cap\partial K_R(tv)_{T+C_{18}+C_{19}}\neq\emptyset}{\Gamma\cap\partial K_R(tv)_{T+C_{18}}\neq\emptyset;\diam(\partial K_R(tv)_{T+C_{18}})\geq1}}\nonumber.
\end{align}
Applying Lemma~\ref{le:q2} conditioned on $\mathcal{F}_{T+C_{18}}$ we obtain that $\prob{\tau^{R}(\Gamma,tv)\leq T}$ can be bounded from above by
\begin{align} 
&\frac{1}{p_1p_2p_3}\prob{\gamma\cap\partial K_R(tv)_{T+C_{18}+C_{19}}\neq\emptyset}=
\frac{1}{p_1p_2p_3}\prob{\gamma_{T+C_{18}+C_{19}}\cap\partial K_R(tv)\neq\emptyset}\nonumber\\
\leq&\frac{1}{p_1p_2p_3}\left( \prob{\tau^R(\gamma,tv)\leq T+C_{18}+C_{19}}+\prob{\diam(\gamma_{T+C_{18}+C_{19}})<1} \right)
\end{align}
where again we used Lemma~\ref{le:reverse}. The fact the we are only considering $T\geq\sqrt{t}$ now shows that for $m\in\mathbb{N}$ there is $\kappa_m^{(11)}\in\R$ such that 
\[
\prob{\diam(\gamma_{T+C_{18}+C_{19}})<1}\vee\prob{\diam(\Gamma_{T+C_{18}})<1}\leq \kappa_m^{(11)} T^{-m}
\]
which completes the proof (choosing $C_{20}:=p_1p_2p_3$ and $C_{17}:=C_{18}+C_{19}$). \hfill$\Box$\\
If we consider now for $t\geq2\frac{C_{17}}{\delta}$ 
\begin{align}
\prob{\tau^R(\Gamma,tv)\leq(1-\delta)t}\leq&\frac{1}{C_{20}}\prob{\tau^R(\gamma^{(t)},tv)\leq(1-\delta)t+C_{17}}+\frac{\kappa_m^{(11)} t^{-m}}{C_{20}}\nonumber\\
\leq&\frac{1}{C_{20}}\prob{\tau^R(\gamma^{(t)},tv)\leq(1-\frac{\delta}{2})t}+\frac{\kappa_m^{(11)} t^{-m}}{C_{20}}\rightarrow0\nonumber
\end{align}
together with $\prob{\tau^R(\Gamma,tv)\leq(1+\delta)t}\geq\prob{\tau^R(\gamma^{(t)},tv)\leq(1+\delta)t}\rightarrow1$
we get that $\tau^{R}(\Gamma,tv)t^{-1}$ converges to $1$ in probability. The diffeomorphic property of the flow of course implies that this convergence holds uniformly in $\gamma\in\mathcal{C}_{R/2}$ if we replace $\Gamma$ by $\gamma$. Corollary~\ref{ko:ggi} also shows that it also holds in $L^p$ for any $p\geq1$. We thus proved the following corollary.
\begin{corollary}\label{ko:konvinwkt} We have   
$\lim_{t\rightarrow\infty} \sup_{\gamma\in\mathcal{C}_R}\prob{\left|\frac{\tau^{R}(\gamma,tv)}{t}-\left\|v\right\|^R\right|>\epsilon}=0$
for $\epsilon>0$ as well as for any $p>0$ that
$\lim_{t\rightarrow\infty}\sup_{\gamma\in\mathcal{C}_R}\expec{\left| \frac{\tau^{R}(\gamma,tv)}{t}-\left\|v\right\|^R \right|^p}=0$.
\end{corollary}
Proof: There is nothing left to show since we ensured that the assertions above do not depend on $R$. \hfill$\Box$\\
We now turn to the last assertion of Theorem~\ref{th:main}.
\begin{lemma}\label{le:osvor} We have for any $\epsilon>0$ that  $\lim_{t\rightarrow\infty}\prob{\mathcal{W}^{R}_{t}(\gamma)\subset (1+\epsilon)t\mathcal{B}}=1$.
\end{lemma}
Proof: For $\eps>0$ we have equivalence with 
\[\lim_{t\rightarrow\infty}\prob{\exists x\in\mathcal{W}^{R}_{t}(\gamma):x\notin(1+\epsilon)t\mathcal{B}}=0.
\]  
So it is enough to show that for $\eps>0$ we have
\[\lim_{t\rightarrow\infty}\prob{\exists x\in\mathbb{R}^{2}:\tau^{R}(\gamma,tx)\leq \frac{t}{1+\epsilon}; \left\|x\right\|^{R}=1}=0.\]
Choose $\delta\ll\epsilon$ and a $\delta$-net on $\partial\mathcal{B}$ denoted by $\left\{v_{j}:j=1,\ldots,N_{\delta} \right\}$. Then we can apply Corollary~\ref{ko:konvinwkt} to obtain
\begin{equation}
\lim_{t\rightarrow\infty}\prob{\exists j\in\{1,\ldots,N_{\delta}\}:\frac{\tau^{R}(\gamma,tv_{j})}{t}<\frac{1}{1+\frac{\epsilon}{2}}}=:\lim_{t\rightarrow\infty}\prob{F_{1}(t)}=0.
\end{equation}
Due to Theorem~\ref{th:cor6} and Lemma~\ref{le:sweep} (similarly to the proof of Theorem~\ref{sa:us}) we have for large $t$ (because for $v\in\mathbb{R}^{2}$ with $\left\|v\right\|^{R}=1$ there is $j$ with $\left\|v-v_{j}\right\|\leq\delta$)
\begin{align}
&\cp{F_{1}(t)}{F_{2}(t)}\nonumber\\
:=&\cp{\exists j\in\{1,\ldots,N_{\delta}\}:\frac{\tau^{R}(\gamma,tv_{j})}{t}<\frac{1}{1+\frac{\epsilon}{2}}} {\exists v\in\partial\mathcal{B}:\frac{\tau^{R}(\gamma,tv)}{t}<\frac{1}{1+\epsilon} }\nonumber\\
\geq&\inf_{v\in\partial\mathcal{B}, \gamma\in\mathcal{C}_{\frac{R}{2}}+tv}\prob{\tilde{\tau}^{\delta t}(\gamma,tv)<\left(\frac{1}{1+\frac{\epsilon}{2}}-\frac{1}{1+\epsilon}\right)t}\nonumber\\
=&\inf_{\gamma\in\mathcal{C}_{\frac{R}{2}}}\prob{\tilde{\tau}^{\delta t}(\gamma,0)<\left(\frac{1}{1+\frac{\epsilon}{2}}-\frac{1}{1+\epsilon}\right)t}\geq1-\kappa_m^{(12)}t^{-m}\label{eq:f1f2}.
\end{align}   
where $\kappa_m^{(12)}\in\mathbb{R}$ exists for $m\in\mathbb{N}$, provided $\delta=\delta(\epsilon)$ is chosen sufficiently small. This implies 
$ 
\limsup_{t\rightarrow\infty}\prob{F_{2}(t)}=0$ and hence Lemma~\ref{le:osvor}.\hfill$\Box$\\ 
The remaining part of Theorem~\ref{th:main} now turns out to be a conclusion from Lemma~\ref{le:osvor} and the fact, that convergence in probability implies a.s. convergence of some subsequence.\hfill $\Box$\\

\end{document}